\documentclass{article}
\usepackage{amsthm,amsfonts,amsmath,amssymb}

\newcommand{\Mohle}{{M\"ohle}}

\newcommand{\ed}{\stackrel{d}{=}}
\newcommand{\prob}{\mathbb P}
\newcommand{\ex}{\mathbb E\,}
\newcommand{\giv}{\,|\,}
\newcommand{\te}{\rightarrow}
\newcommand{\eq}{\begin{equation}}
\newcommand{\en}{\end{equation}}
\newcommand{\FM}{{\rm FM}}

\newcommand{\re}[1]{\mbox{(\ref{#1})}}
\newcommand{\rem}[1]{\mbox{\rm (\ref{#1})}}

\newcommand{\Nat}{\Bbb N}

\def\endpf{\hfill $\Box$ \vskip0.5cm}
\def \proof{\noindent{\it Proof.\ }}
\def \Lev{L\'evy}
\theoremstyle{plain}

\newtheorem{theorem}{\large Theorem}
\newtheorem{proposition}[theorem] {\large Proposition}
\newtheorem{definition}[theorem]{\large Definition}

\newtheorem{lemma}[theorem]{\large Lemma}

\begin{document}

\title{Exchangeable partitions derived from Markovian coalescents
\thanks{Research supported in part by N.S.F. Grant DMS-0405779}
}
\author{Rui Dong
\thanks{University of California, Berkeley; e-mail ruidong@stat.Berkeley.EDU}  
\hspace{.2cm}
and 
Alexander Gnedin
\thanks{Utrecht University; e-mail gnedin@math.uu.nl}
\hspace{.2cm}
and
Jim  Pitman\thanks{University of California, Berkeley; e-mail pitman@stat.Berkeley.EDU}  
\\
}

\date{\today}
\maketitle

\begin{abstract}
\noindent
Kingman derived the Ewens sampling formula for random partitions
describing the genetic variation in a neutral mutation model
defined by a Poisson process of mutations along lines of descent
governed by a simple coalescent process, and observed that similar
methods could be applied to more complex models.
\Mohle\ described the recursion which determines 
the generalization of the Ewens sampling formula in the situation when
the  lines of descent are governed by a $\Lambda$-coalescent,
which allows multiple mergers.
Here we show that the basic integral representation of 
transition rates for the $\Lambda$-coalescent is forced by 
sampling consistency under
more general assumptions on the coalescent process.  Exploiting an 
analogy with the theory of regenerative  partition structures, we provide
various characterizations of the associated partition structures in 
terms of discrete-time Markov chains.
\end{abstract}

\section{Introduction} 
The theory of random coalescent processes starts from Kingman's 
series of papers \cite{kingman82a,  kingman82c, kingman82b} in 1982. The idea comes from biological studies
 for genealogy of haploid model \cite{canning74}: given a large population with many generations, you track backward in time 
the family history of each individual in the current generation. As you track further, 
the family lines coalesce with each other, eventually all terminating at a common ancestor of current generation.  The same mathematical process may be 
interpreted in other way as describing collisions of an aggregating system of physical particles.
In Kingman's coalescent process \cite{kingman82a}, each collision only involves two parts. 
This idea is extended to coalescent with multiple collisions in \cite{lambda, sagitov}, where every collision can involve two or more parts. 
This model is further developed into the theory of coalescent with simultaneous multiple collisions in \cite{jason00, moehlesagitove01}. 
See 
\cite{tavare84,watter84,evanspitman98, boltszni98, sagitov03, bertoingold04, donggoldmartin05} for related developments.

Kingman \cite{kingman82b}  indicated a basic connection between random partitions of natural interest in genetics, and coalescent processes.
Suppose in the haploid case the family line of current generation is modeled by Kingman's coalescent, and the mutations are applied along the family lines by using a Poisson process with rate $\theta/2$ for some non-negative 
 number $\theta$. Define a partition by saying that two individuals are in the 
same block if there is no mutation along their family lines before they coalesce. Then the resulting random partition is governed by the Ewens sampling formula 
with parameter $\theta$. See   \cite[Section 5.1, Exercise 2]{CSP} and \cite{BertoinBook, Nordborg} for review and more on this idea. 
Recently, M{\"o}hle \cite{moehle1} applied 
this idea to the genealogy tree modeled by coalescents with multiple collisions and simultaneous multiple collisions. 
He studied the resulting family of 
partitions, and derived a recursion which determines them. In \cite{moehle2}, M{\"o}hle showed that the partition derived 
from 
coalescent with multiple collisions is regenerative in the sense of \cite{RCS, RPS} if and only if the underlying coalescent 
is Kingman's coalescent or 
a hook case, corresponding to the extreme cases when the 
characterization measure $\Lambda$ of coalescent with multiple collisions concentrates at $0$ or $1$, 
respectively. 
In particular, 
the intersection of \Mohle's family of partitions with Pitman's two-parameter family 
is the one-parameter Ewens' family. 

Here we offer a different approach to the family of random partitions
generated by Poisson marking along the lines of descent of 
a $\Lambda$-coalescent.
We study partitions with an additional feature, assigning each part one of two possible states: active or frozen.
We introduce a new class of continuous time  partition-valued coalescent 
processes, called {\em coalescents with freeze}, 
which are characterized by an underlying  measure determining collision
rates, together with a freezing rate. Every coalescent with freeze has a terminal state with all blocks frozen, 
called the {\em final partition} of this process, whose distribution
is characterized by the recursion of M{\"o}hle \cite{moehle1}.
In the spirit of \cite{RCS, RPS}, we focus here on the discrete time chains 
embedded in the coalescent with freeze, and from the consistency of their transition operators 
we derive a backward recursion satisfied by the decrement 
matrix, analogous
to \cite[Theorem 3.3]{RCS}. This decrement matrix determines the partition through M\"ohle's recursion. As in \cite{RCS}, 
we use algebraic methods to derive an integral representation for the decrement matrix.  Also, adapting an  idea from \cite{RPS}, we 
establish a uniqueness result by constructing another Markov chain, with state
space the set of partitions of a finite set, whose unique stationary 
distribution is the law of the final partition restricted to this set.
We analyze in detail the case of coalescent with freeze when no simultaneous multiple collisions 
are permitted, leaving the more general case to another paper.

\par The remaining part of the paper is organized as following. Some notations and background are introduced in Section \ref{2}, 
together with a review of M\"ohle's result. 
In Section \ref{3} the coalescent with freeze is defined and the relation between our method and M\"ohle's method is discussed. In Section \ref{4} we 
detail the study of coalescent with freeze in terms of the freeze-and-merge (FM) operators of the embedded finite discrete chain, 
whose consistency with sampling derives a backward recursion for the decrement matrix. In Section 5 the Markov chain with 
sample-and-add (SA) operation is introduced, and the law of 
the partition in our study is identified as the unique stationary distribution 
of this chain. 
In Section 6 we derive the integral representation for an infinite decrement matrix.  This gives another approach to \Mohle's partitions via consistent 
freeze-and-merge chains, which may be seen as  discrete-time jumping 
processes associated with the $\Lambda$-coalescent with freeze.
Section \ref{pos} provides an alternate approach to the representation of
an infinite decrement matrix in terms of a positivity condition on a single 
sequence.
Section \ref{7} offers some results about the structure of the random
set of freezing times derived from a coalescent with freeze. Finally,
in Section \ref{8} we point out some striking parallels with our previous 
work on regenerative partition structures, which guided this study.
Section \ref{9} mentions briefly some further parallels with 
the theory of homogenous and self-similar Markovian fragmentation 
processes due to Bertoin \cite{BertoinBook}.

\section{Some notation and background}\label{2}

Following the notations of \cite{CSP}, for any finite set $F$, a partition of $F$ into $\ell$ blocks, also called a {\em finite set partition},
is an unordered collection of non-empty disjoint sets $\{A_1,\ldots,A_\ell\}$ whose union is $F$. In particular we consider partitions 
of the set $[n]:=\{1,2,\ldots,n\}$ for $n\in \mathbb{N}$. We use $\mathcal{P}_{[n]}$ to denote the set of all partitions of $[n]$.
A {\em composition} of the positive integer $n$ is an ordered sequence of positive integers $(n_1,n_2,\ldots,n_\ell)$ with $\sum_{i=1}^\ell n_i=n$, where $\ell\in \mathbb{N}$ is number of parts. We use $\mathcal{C}_n$ to denote the set of all compositions of $n$, and
$\mathcal{P}_n$ to denote the set of non-increasing compositions of $n$, also called {\em partitions of $n$}.

Let $\pi_n=\{A_1,A_2,\ldots,A_\ell\}$ denote a generic partition of $[n]$; we may write $\pi_n\vdash [n]$ to indicate this fact.
The {\em shape function} from partitions of the set $[n]$ to partitions of the positive integer $n$ is defined by
\eq
\mathtt{shape}(\pi_n)=(|A_1|,|A_2|,\ldots,|A_\ell|)^{\downarrow}
\en 
where $|A_i|$ is the size of block $A_i$ which represents the number of elements in 
the block,
 and ``$\downarrow$'' means arranging the sequence of sizes in non-increasing order.
 
A random partition $\Pi_n$ of $[n]$ is a random variable taking values in $\mathcal{P}_{[n]}$. It is called exchangeable 
if its distribution is invariant under the action on partitions of $[n]$ by 
the symmetric group of 
permutations of $[n]$. 
Equivalently, the distribution of $\Pi_n$ is given by the formula
\eq
\mathbb{P}(\Pi_n=\{A_1,A_2,\ldots,A_\ell\}) = p_n(|A_1|, |A_2|,\ldots,|A_\ell|)
\en
for some symmetric function $p_n$ of compositions of $n$.
 We call $p_n$ the {\em exchangeable partition probability function (EPPF)} of $\Pi_n$.

An {\em exchangeable random partition of $\Nat$} is
a sequence of exchangeable set partitions $\Pi_{\infty}=(\Pi_n)_{n=1}^\infty$ with $\Pi_{n}\vdash [n]$, subject to the
{\em consistency condition}
\eq\label{infinitepartition}
\Pi_n |_{m} = \Pi_m,
\en
where the restriction operator $|_{m}$ acts on $\mathcal{P}_{[n]}$, $n>m$, by deleting elements $m+1,m+2,\ldots,n$. 
The distribution of such an exchangeable random partition of $\Nat$ is determined by the
function $p$ defined on the set of all integer compositions $\mathcal{C}_\infty:=\cup_{i=1}^\infty \mathcal{C}_i$, which coincides with the 
EPPF $p_n$ of $\Pi_n$ when acting on $\mathcal{C}_n$. This function $p$ is called the infinite EPPF associated
with $\Pi_{\infty}=(\Pi_n)_{n=1}^\infty$.
The consistency condition \re{infinitepartition} translates into the following {\em addition rule} for the EPPF $p$: for each positive integer 
$n$ and each composition $(n_1,n_2,\ldots,n_\ell)$ of $n$,
\eq\label{sps}
p(n_1,n_2, \ldots,n_\ell) = p(n_1,n_2, \ldots, n_\ell, 1)+\sum_{i=1}^\ell p(n_1, \ldots, n_i+1, \ldots, n_\ell)
\en
where $(n_1, \ldots, n_i+1, \ldots, n_\ell)$ is formed from $(n_1, \ldots, n_\ell)$ by adding $1$ to $n_i$. 
Conversely, if a nonnegative function $p$ on compositions satisfies  (\ref{sps}) and the normalization condition $p(1)=1$,
then by Kolmogorov's extension theorem there exists an exchangeable random  partition $\Pi_\infty$ with EPPF $p$.

\par
Similar definitions apply to a finite sequence of consistent exchangeable random set partitions 
$(\Pi_m)_{m=1}^n$ with $\Pi_{m}\vdash [m]$, where $n$ is some fixed positive integer. The {\it finite} EPPF $p$ of such a sequence can be defined 
as the unique recursive extension of $p_n$ 
by the 
 addition rule \re{sps} to all compositions
$(n_1,n_2,\ldots,n_\ell)$ of $m < n$.

\par
Let $\mathcal{P}_{\infty}$ be the set of all partitions of $\Nat$. 
We identify each $\pi_\infty \in \mathcal{P}_{\infty}$ as the sequence 
$(\pi_{1},\pi_{2},\ldots)\in \mathcal{P}_{[1]}\times \mathcal{P}_{[2]}\times\cdots$, where $\pi_{n}=\pi_\infty|_n$ is the restriction 
of $\pi_{\infty}$ to $[n]$ by deleting all elements bigger than $n$. Give $\mathcal{P}_{\infty}$ the topology it inherits as a subset 
of $\mathcal{P}_{[1]}\times \mathcal{P}_{[2]} \times\cdots$ with the product of discrete topologies, so the space $\mathcal{P}_{\infty}$ is compact 
and metrizable. 
Following \cite{evanspitman98, kingman82a,lambda}, 
call a $\mathcal{P}_{\infty}$-valued stochastic process $(\Pi_{\infty}(t), t\ge 0)$ a {\it coalescent} if it has c\`adl\`ag paths 
and $\Pi_{\infty}(s)$ is a {\it refinement} of $\Pi_{\infty}(t)$ for every $s<t$. 
For a non-negative finite measure $\Lambda$ on the Borel subsets of $[0,1]$,
a {\em $\Lambda$-coalescent} is a $\mathcal{P}_{\infty}$-valued Markov coalescent $(\Pi_{\infty}(t), t\ge 0)$ whose restriction 
$(\Pi_{n}(t), t\ge 0)$ to $[n]$ is for each $n$ a Markov chain such that when $\Pi_{n}(t)$ has $b$ blocks, 
each $k$-tuple of blocks of $\Pi_{n}(t)$ is merging to form a single block at rate $\lambda_{b,k}$, where
\eq\label{lambdaint}
\lambda_{b,k}=\int_0^1x^{k-2}(1-x)^{b-k}\Lambda(dx)\ \ \ \ \ \ (2\le k\le b<\infty).
\en
The measure $\Lambda$ which characterizes the coalescent is derived from the consistency requirement, 
that is for any positive integers $0<m<n<\infty$, and $\pi_n\vdash [n]$, the restricted process $(\Pi_n(t)|_{m},\,t\ge 0)$ given $\Pi_{n}(0)=\pi_{n}$ 
has the same law as $(\Pi_{m}(t), t\ge 0)$ given $\Pi_{m}(0)=\pi_{n}|_{m}$. This condition is fulfilled if and only if  
the array of rates $(\lambda_{b,k})$ satisfies
\eq\label{pitmanlambda}
\lambda_{b,k}=\lambda_{b+1,k}+\lambda_{b+1,k+1}~~~~~~~~~~~~~~~~(2\le k\le b<\infty).
\en
The integral representation (\ref{lambdaint}) can be derived from \re{pitmanlambda} via 
de Finetti's theorem \cite[Lemma 18]{lambda}. 

\par
When $\Lambda=\delta_0$, this reduces to Kingman's coalescent \cite{kingman82a, kingman82b, kingman82c} with only binary merges. 
When $\Lambda$ is the uniform distribution on $[0,1]$, the coalescent is the Bolthausen-Sznitman coalescent \cite{boltszni98}. 
In \cite{jason00} this construction is further developed to build the $\Xi$-coalescent where the measure $\Xi$ on infinite simplex characterizes 
the rates of simultaneous multiple collisions. 

\par
M{\"o}hle \cite{moehle1} studied the following generalization of Kingman's model \cite{kingman82b}. Take a genetic sample of 
$n$ individuals
from a large population and label them as $\{1,2,\ldots,n\}$. Suppose the ancestral lines of these $n$ individuals evolve by the rules of
 a $\Lambda$-coalescent, and that given the genealogical tree, whose branches are the ancestral lines of these individuals, 
mutations occur along the ancestral lines according to a Poisson point process with rate $\rho>0$.
The infinite-many-alleles model is assumed, which means that when a gene mutates,  a brand new type appears.
Define a random partition of $[n]$ by declaring individuals $i$ and $j$ to be in the same block if and only if they are of the same type,
that is either $i=j$ or there are no mutations along the ancestral lines of $i$ and $j$ before these lines coalesce. 
These random partitions are exchangeable, and consistent as $n$ varies.
The EPPF of this random partition is the unique solution $p$ with $p(1) = 1$ of 
{\em \Mohle's recursion}:
for each  positive integer $n$ and each composition $(n_1,n_2,\ldots,n_\ell)$ of $n$,
\eq\label{rec1}
p(n_1,n_2,\ldots,n_\ell)=
{q(n:1)\over n}\sum_{j:n_j=1} p(\ldots,\widehat{n_j},\ldots)+ 
\sum_{k=2}^n q(n:k)\sum_{j: n_j\ge k} \frac{{n_j\choose k}}
{{n\choose k}}p(\ldots, n_j-k+1,\ldots),
\en
where $(\ldots,\widehat{n_j},\ldots)$  
is formed from $(n_1,n_2,\ldots,n_\ell)$ by removing part $n_j$, $(\ldots, n_j-k+1,\ldots)$ is formed from  $(n_1,n_2,\ldots,n_\ell)$ 
by only changing $n_j$ to $n_j-k+1$, and 
$q(b:k)$ 
is the stochastic  matrix
\eq\label{M}
 q(b:k)={\Phi(b:k)\over \Phi(b)}\,~~~~~(1\leq k\leq b \leq n),
\en
where
\begin{eqnarray}
\Phi(b:1)&=&\rho b\,,\label{Phi1}\\
\Phi(b:k)&=&{b \choose k } \lambda_{b,k} = {b\choose k}\int_0^1 x^{k-2}(1-x)^{b-k}\Lambda({\rm d}x)\,~~~~~(2\leq k\leq b)\,,\label{Phi2}\\
\Phi(b)&=&\sum_{k=1}^b\Phi(b:k)=\int_0^1{1-(1-x)^b-bx(1-x)^{b-1}\over x^2}\,\,\Lambda({\rm d}x)\,+\rho b
\label{Phi3}.
\end{eqnarray}
If at some time $t \ge 0$  there are exactly $b$ lines of descent 
whose associated genealogical trees of depth $t$ contain no mutations,
then $\Phi(b:1)$ is the total rate of mutations along one of these $b$ lines, 
$\Phi(b:k)$ is the total rate of $k$-fold merges among these lines,
and $\Phi(b)$ is the total rate of events of either kind.

\par
M{\"o}hle \cite{moehle1} derived the recursion 
\re{rec1}
by conditioning on whether the first event met tracing 
back in time from the current generation is a mutation or collision. On the left side of \re{rec1},
$p(n_1,n_2,\ldots,n_\ell)$ is the probability of ending up with any {\em particular} partition $\pi_n$ of the set $[n]$ into
$\ell$ blocks of sizes $(n_1,n_2,\ldots,n_\ell)$.
On the right side, $q(n:1)$
is the chance that starting from the current generation, one of the $n$ genes mutates before any collision; for this to happen together with the
specified partition of $[n]$, the individual with this gene must be chosen from those among the singletons of $\pi_n$, with chance  $1/n$ for 
each different choice, and after that the restriction of the coalescent process to a subset of $[n]$ of size $n-1$ must end up generating the 
restriction of $\pi_n$ to that set.
Similarly, $q(n:k)$ is the chance that the first event met is $k$ out of $n$ genes coalescing to the same block.
Again, the $k$ individuals bearing these $k$ genes must be chosen from a block of $\pi_n$ of size $n_j \ge k$, 
so the chance for possible choices from a block with size $n_j$ is ${n_j \choose k}/{n\choose k}$,
and given exactly which $k$ individuals are chosen, the restriction of the coalescent process to some set of $n-k+1$ lines of descent
must end up generating  a particular partition of these $n-k+1$ lines into sets of sizes $(\ldots, n_j-k+1,\ldots)$.
The multiplication of various probabilities is justified by the strong Markov property of the $\Lambda$-coalescent at the time of the first event,
and by the special symmetry property that lines of descent representing blocks of individuals coalesce according to the same dynamics as if they
were singletons.

\par In this paper we step back from these detailed dynamics of the $\Lambda$-coalescent with mutations to consider the following questions 
related to M{\"o}hle's recursion \re{rec1} and associated partition-valued processes. 
We choose to ignore the special form \re{M} of the matrix 
$(q(n:k);~ 1 \le k \le  n<\infty)$ derived from the $(\Lambda, \rho)$, 
and analyse M{\"o}hle's recursion \re{rec1} as an abstract relation between a stochastic
 matrix $q$ and a function of compositions
$p$. In particular, we ask the following questions:
\begin{enumerate}
\item For which probability distributions $q(n:k),\,1\leq k\leq n$, on $[n]$
 is \Mohle's recursion \re{rec1} satisfied by the EPPF  $p$ of
some exchangeable random partition of $[n]$, and is this $p$ uniquely determined?
\item How can such random partitions be characterized probabilistically?
\item Can such random partitions of $[n]$ be consistent as $n$ varies for any other $q$ besides $q$ derived from $(\Lambda, \rho)$ as above?
\end{enumerate}
We stress that in the first two questions the recursion \re{rec1} is only required to hold for 
a single value of $n$, while in the third question \re{rec1} must hold for all $n=1,2,\ldots$.
The answer to the first question is that for each fixed probability distribution $q(n:k),\, 1 \le k \le  n$, on  $[n]$,
\Mohle's recursion \re{rec1} determines a unique EPPF $p$ for an exchangeable random partition of $[n]$ (Theorem \ref{thatsit}).
Answering the second question, we characterize the distribution of this random partition in two different ways: firstly as the terminal state of a
discrete-time Markovian coalescent process, the {\em freeze-and-merge chain} introduced in Section \ref{4}, and
secondly as the stationary distribution of a partition-valued Markov chain with quite a different transition mechanism, the
{\em sample-and-add chain} introduced in Section \ref{5}.
The answer to the third question is positive if we restrict $n$ to some bounded range of values, for some but not all $q$
(see Section \ref{4}), but negative if we require consistency for all $n$ (Theorem \ref{MAIN}): 
if an infinite EPPF $p$ solves \Mohle's recursion \re{rec1} for all
$n$ for some triangular matrix $q$ with non-negative entries, then $q$ must have the form \re{M}
for some $(\Lambda,\rho)$.
\par
We were guided in this analysis by a remarkable parallel 
between this theory of finite and infinite partitions subject to \Mohle's recursion \re{rec1} and 
the theory of {\em regenerative partitions} developed in \cite{RCS, RPS}. 
Following the terminology in \cite{RCS,RPS}, we call a triangular stochastic matrix 
a {\em decrement matrix}. We use the notation 
$q_n=(q(b:k); ~\,1\leq k\leq b\leq n)$
or $q_\infty=(q(n:k);~ \,1\leq k\leq n<\infty)$
to indicate whether we wish to consider finite or infinite matrices.
Thus, the entries of a decrement matrix are nonnegative and  
satisfy
 $\sum_{k = 1}^b q(b:k) = 1$ for all $b$ in the required range.
In present notation, the characteristic property of a regenerative partition  is that its EPPF $p$ satisfies
\eq
\label{recregen}
p(n_1,n_2,\ldots,n_\ell)=
\sum_{j=1}^\ell \frac{1}{ {n \choose n_j } }\, q(n:n_j)\, p(\ldots,\widehat{n_j},\ldots)
\en
for some decrement matrix $q=q_\infty$.  
The main results of \cite{RCS,RPS} gave similar answers to the above questions for this recursion instead 
of \Mohle's recursion \re{rec1}.

\par There is an important distinction between the recursion 
(\ref{sps}) on the one hand and (\ref{rec1}) and \re{recregen} on the other hand.
The recursion (\ref{sps}) has many solutions since it
is a {\it backward} recursion, from larger values of $n$ to smaller.
By contrast, 
both (\ref{rec1}) and (\ref{recregen}) are {\it forward} recursions, from smaller values of $n$ to larger.
Consequently it is obvious that given an arbitrary infinite decrement matrix $q_\infty$, 
each of the recursions {\rm (\ref{rec1})} and \re{recregen} has a unique solution $p$ with the initial value $p(1)=1$.
Moreover, it is clear that each of these functions $p$ can be written as a linear combination of products of entries
of the $q_\infty$ matrix.
\par
To illustrate the close parallel between the two recursions \re{rec1} and \re{recregen},
we list the first few values of $p$ in terms of the decrement matrix $q$, first for \Mohle's recursion (\ref{rec1}):

\begin{align*}
p(1)&=  1,\\
p(2)&=  q(2:2),\\
p(1,1)&= q(2:1),\\
p(3)&= q(3:3)+q(3:2)q(2:2),\\
p(2,1)&= p(1,2)\\
&= {1\over 3}q(3:2)q(2:1)+{1\over 3}q(3:1)q(2:2),\\
p(1,1,1)&= q(3:1)q(2:1),\\
p(4)&= q(4:4)+q(4:3)q(2:2)+q(4:2)q(3:3)+q(4:2)q(3:2)q(2:2),
\end{align*}
\begin{align*}
p(3,1)&= p(1,3)\\
&= \frac{1}{4}q(4:3)q(2:1)+\frac{1}{6}q(4:2)q(3:2)q(2:1)+\frac{1}{2}q(4:2)q(3:1)q(2:2)+\frac{1}{4}q(4:1)q(3:3)\\
&\ \ \ +\frac{1}{12}q(4:1)q(3:2)q(2:2),\\
p(2,1,1)&= p(1,2,1)=p(1,1,2)\\
&= \frac{1}{6}q(4:2)q(3:1)q(2:1)+\frac{1}{6}q(4:1)q(3:2)q(2:1)+\frac{1}{6}q(4:1)q(3:1)q(2:2),\\
p(1,1,1,1)&= q(4:1)q(3:1)q(2:1).
 \end{align*}
Note that for a general transition matrix $q$ these functions $p$  may not  be consistent as $n$ varies, meaning
that \re{sps} may fail. A condition on $q_\infty$ equivalent to consistency of $p$ will be described later in Lemma \ref{l5}.

Similarly, the first few values of the $p$ determined by a decrement matrix $q$ via the recursion \re{recregen}
associated with a regenerative partition structure are:
\begin{align*}
p(1)&=1,\\
p(2)&=q(2:2),\\
p(1,1)&=q(2:1),\\
p(3)&=q(3:3),\\
p(2,1)&=p(1,2)\\
&=\frac{1}{3}q(3:2)+\frac{1}{3}q(3:1)q(2:2),\\
p(1,1,1)&=q(3:1)q(2:1),\\
p(4)&=q(4:4),\\
p(3,1)&=p(1,3)\\
&=\frac{1}{4}q(4:3)+\frac{1}{4}q(4:1)q(3:3),\\
p(2,1,1)&=p(1,2,1)=p(1,1,2)\\
&=\frac{1}{6}q(4:2)q(2:1)+\frac{1}{6}q(4:1)q(3:2)+\frac{1}{6}q(4:1)q(3:1)q(2:2),\\
p(1,1,1,1)&=q(4:1)q(3:1)q(2:1).
\end{align*}
Looking at these displays, both similarities and differences may be observed.
In particular, the formulas for singleton partitions $(1,1,\ldots,1)$ are identical. 
As is to be expected, the simpler recursion \re{recregen} for regenerative partitions generates simpler algebraic expressions
than  \Mohle's recursion \re{rec1}. See \cite[Equation (16)]{RPS} (reproduced as (\ref{EPFRegen}) below) for the general formula for the shape 
function associated with \re{recregen}.

\par In principle, the recursions \re{rec1} and \re{recregen} have probabilistic meaning for arbitrary decrement matrix $q$,
since they determine a sequence of exchangeable partitions of $[n]$'s for $n$ in some finite or the infinite range.
Distributions of these partitions are obtained algebraically as above, by fully expanding $p$ through $q$. 
However, typically these partitions of $n$ are not consistent with respect to restrictions, so in the infinite case 
they might not determine the distribution of a partition of $\Nat$.

\section{Coalescents with freeze} \label{3}

To provide a natural generalization of partition structures derived from a coalescent with Poisson mutations along the branches of a genealogical tree,
we consider the structure of a partition of a set (respectively, of an integer) with each of its blocks (or parts) assigned one of two possible conditions, 
which we call 
{\it active} and {\it frozen}.  We call such a combinatorial object a {\it partially frozen partition}  of a set or of an integer, as the case may be.
Ignoring the conditions of the blocks
of a partially frozen partition  $\pi^*$
{\em induces} an ordinary partition $\pi$.
As special cases of partially frozen partititions, we include the possibilty that all blocks may be active, or all frozen.
We use the symbol $\Sigma_n^*$ for the pure singleton partition of $[n]$ with all blocks active, and $\Sigma_\infty^*$ for the 
sequence $(\Sigma_n^*)_{n=1}^\infty$. 
The {\it *-shape} 
of a partially frozen partition $\pi_n^*$ of $[n]$ is 
the corresponding partially frozen partition of $n$, and the ordinary
shape is defined in terms of the  induced partition $\pi_n$.

For each positive integer $n$, we denote  $\mathcal{P}_{[n]}^*$ the set of all partially frozen 
partitions of $[n]$. 
Let $\mathcal{P}_{\infty}^*$ be the set of all partially frozen 
partitions of $\Nat$. 
We identify each element $\pi_\infty^*\in \mathcal{P}_{\infty}^*$ as the sequence 
$(\pi_1^*,\pi_2^*,\ldots)\in \mathcal{P}_{[1]}^*\times \mathcal{P}_{[2]}^*\times\cdots$, where $\pi_n^*$ is $\pi_\infty^*|_{n}$ 
the restriction of $\pi_\infty^*$ to $[n]$. Endowing $\mathcal{P}_{\infty}^*$ with the topology it inherits as a subset of 
$\mathcal{P}_{[1]}^*\times \mathcal{P}_{[2]}^*\times\cdots$, the space $\mathcal{P}_{\infty}^*$ is compact and metrizable. 
We call a random partially frozen partition of $[n]$ exchangeable
if its
 distribution is invariant under the action of permutations of $[n]$. 
Similarly to \cite{evanspitman98,kingman82a}, call a $\mathcal{P}_{\infty}^*$-valued stochastic process 
$(\Pi_{\infty}^*(t), t\ge 0)$ a {\it coalescent} if it has c\`adl\`ag paths and $\Pi_{\infty}^*(s)$ is a {\it *-refinement} of $\Pi_{\infty}^*(t)$
 for every $s<t$, 
 meaning that
the induced partition $\Pi_{\infty}(s)$ is a refinement of $\Pi_{\infty}(t)$ and the set of frozen blocks of $\Pi_{\infty}^*(s)$ is 
a subset of the set of frozen blocks 
of $\Pi_{\infty}^*(t)$.

\par 
The construction of an exchangeable random partition of  $\Nat$ by cutting branches of the merger-history tree of a $\Lambda$-coalescent
$(\Pi_{\infty}(t),t \ge 0)$ 
by mutations with rate $\rho$ can now be formalized as follows.
For each $i \in \Nat$ let $\tau_i$ denote the random time at which a mutation first occurs along the line of descent to leaf $i$ of the tree,
and declare the block of $\Pi_{\infty}(t)$ containing $i$ to be active if $\tau_i > t$ and frozen if $\tau_i \le t$.
This defines a $\mathcal{P}_{\infty}^*$-valued Markov process $(\Pi_{\infty}^*(t), t \ge 0)$.
As $t \te \infty$ the state $\Pi_{\infty}^*(t)$ approaches a limit $\Pi_{\infty}^*(\infty)$
with all blocks frozen.
This is the
exchangeable random partition generated by the exchangeable sequence of random variables $(\tau_i, i \in \Nat)$,
meaning that two integers $i$ and $j$ are in the same block 
of $\Pi_{\infty}^*(\infty)$
iff $\tau_i=\tau_j$.
Assuming that $\Pi_{\infty}^*(0)=\Sigma^*_\infty$, it should be clear that the EPPF of $\Pi_{\infty}^*(\infty)$ is that
defined by \Mohle's recursion \re{rec1}.
The following two theorems present  more formal statements.

\begin{theorem}\label{coalescentwithfreeze}
Let $(\lambda_{b,k},2\le k\le b<\infty)$, $(\rho_n, 1\le n<\infty)$ be two arrays of non-negative real numbers. There exists for 
each $\pi_\infty^*\in \mathcal{P}_{\infty}^*$ a $\mathcal{P}_{\infty}^*$-valued coalescent $(\Pi_{\infty}^*(t), t\ge 0)$ with 
$\Pi_{\infty}^*(0)=\pi^*_\infty$, for each $n$ whose restriction $(\Pi_{n}^*(t), t\ge 0)$ to $[n]$ is a $\mathcal{P}_{[n]}^*$-valued Markov chain starting from 
$\pi_n^*=\pi_\infty^*|_{n}$, and evolving with the rules:
\begin{itemize}
\item at each time $t \ge 0$, conditionally given $\Pi_{n}^*(t)$ with $b$ active blocks, each $k$-tuple of active blocks of $\Pi_{n}^*(t)$ is 
merging to form a single active block at rate $\lambda_{b,k}$, and 

\item each active block turns into a frozen block at rate $\rho_{n,b}$,
\end{itemize}
if and only if  the integral representation \rem{lambdaint} holds
for some non-negative  finite measure $\Lambda$ on the Borel subsets of $[0,1]$, and
$\rho_{n,b}=\rho$
for some non-negative real number $\rho$. 
This $\mathcal{P}_{\infty}^*$-valued process $(\Pi_\infty^*(t), t\ge 0)$ directed by
 $(\Lambda,\rho)$ is a strong Markov process. 
For $\rho=0$, this process reduces to the $\Lambda$-coalescent, and for $\rho >0 $ the process is obtained by superposing Poisson marks
at rate $\rho$ on the merger-history tree of a $\Lambda$-coalescent, and freezing the block containing $i$ at the time of the first mark 
along the line of descent of   $i$ in the merger-history tree.
\end{theorem}
\proof
Just as in \cite{lambda}, consistency of the rate descriptions for different $n$
implies that \re{pitmanlambda} holds, hence the integral representation \re{lambdaint}, and equality of the
$\rho_{n,b}$'s is also obvious by consistency.
\endpf
\begin{definition}\label{def1}
{\rm 
Call this $\mathcal{P}_{\infty}^*$-valued Markov process directed  by a non-negative integer $\rho$ and a non-negative finite measure $\Lambda$ on $[0,1]$ 
the {\it $\Lambda$-coalescent freezing at rate $\rho$},
or the {\it $(\Lambda,\rho)$-coalescent} for short.
 Call a $(\Lambda,\rho)$-coalescent starting from state $\Sigma_\infty^*$ a 
{\it standard $\Lambda$-coalescent freezing at rate $\rho$}, where $\Sigma_\infty^*$ is the pure singleton partition with all blocks active.
}
\end{definition}

Consider the finite coalescent with freeze $(\Pi^*_n(t), t\ge 0)$ which is  the restriction of a standard 
$\Lambda$-coalescent freezing at rate $\rho$ to $[n]$. According to the description above, 
all active blocks will coalesce by the rules of a  $\Lambda$-coalescent, except that every active block enters the frozen condition
at rate $\rho$, 
and after that the block  will stay frozen forever. Hence it is clear that as long as the freezing rate $\rho$ is positive, in finite time the process
$(\Pi_n^*(t),t\ge 0)$ 
will eventually reach a
 {\it final partition} $E_n^*$, with all of its blocks in the  frozen condition.

\par 
Now recall M\"ohle's model \cite{moehle1} as reviewed in Section \ref{2}. The ancestral lines of $n$ labeled genes of current generation 
coalesce as a $\Lambda$-coalescent, and mutations happen along each ancestral line as Poisson point process with rate $\rho>0$.
Hence the final
partition of $[n]$ is defined so that if 
the ancestral line of 
an individual is interrupted  by a mutation before the line coalesces with any other ancestral lines,  
the individual will be a singleton in the partition. This corresponds to the idea of freezing here: 
tracing evolution of a particle  starting from time $0$,  
if a particle freezes before coalescing with others, it will enter as a singleton block in the final partition of the process. 

To detail the study, let us look at the discrete chain embedded in $\Lambda$-coalescent freezing at rate $\rho$. 
By the definition, for each time $t\ge 0$, $\Pi_n^*(t)$ is a partially frozen exchangeable random partition of $[n]$, hence 
its induced form $\Pi_n(t)$ gives an exchangeable random  partition of $[n]$. So does the final partition $E_n^*=\Pi_{n}^*(\infty)$ 
and its induced 
form $E_n$. Set $E_\infty^*:=(E_n^*)$ as the final partition of $(\Pi_\infty^*(t), t\ge 0)$, and denote its induced partition as $E_\infty=(E_n)$. 
The following facts can be read from the existence of $(\Pi_\infty^*(t),t\ge 0)$ and \Mohle's analysis recalled around \re{rec1}.

\begin{theorem} {\em ( M\"ohle \cite[Theorem 3.1]{moehle1})}
The induced final partition $E_\infty=(E_n)_{n=1}^\infty$ of a standard $\Lambda$-coalescent freezing at rate $\rho>0$ is an exchangeable infinite 
random partition of $\Nat$ whose
EPPF $p$ is the unique solution of \Mohle's recursion \rem{rec1} with coefficients from 
the infinite decrement matrix $q_\infty$ defined through $(\Lambda,\rho)$ as in \rem{M}.
\end{theorem}

\section{Freeze-and-merge operations} \label{4}
Given a stochastic process $X$ indexed 
by a continous time parameter $ t \ge 0$, assuming  $X$ has right continuous piecewise constant paths,
the {\em jumping process derived from $X$} is the discrete-time process
$$
\widehat{X}=(\widehat{X}(0),\widehat{X}(1),\ldots) = ( X(T_0), X(T_1), X(T_2), \ldots )
$$
where $T_0 := 0$ and $T_k$ for $k \ge 1$ is the least $t > T_{k-1}$ such that $X(t) \neq X(T_{k-1})$,
if there is such a $t$, and $T_k = T_{k-1}$ otherwise.
The processes $X$ of interest here will ultimately arrive in some absorbing state, and then so too
will $\widehat{X}$.
In particular, the finite coalescent with freeze   $(\Pi^*_n(t), t\ge 0)$, obtained by
restriction to $[n]$ of a $\Lambda$-coalescent freezing at positive rate $\rho$, 
is a Markov chain with transition rate 
${b\choose k} \lambda_{b,k}$ for a $k$-merge and rate $b\rho$ for a freeze,
where $b$ is the number of active blocks at time $t$ and 
the $\lambda_{b,k}$'s are as in 
\re{lambdaint}; while
the jumping process
$\widehat{\Pi}^*_n$ is then a Markov chain governed by the following
{\it freeze-and-merge} operation  ${\rm FM}_n$, which acts on a generic partially frozen partition 
$\pi^*_n$ of $[n]$ as follows:  if $\pi^*_n$ has $b>1$ active blocks then
\begin{itemize}
\item  with probability $q(b:k)$ some $k$ of $b$ active blocks are chosen uniformly at random and merged into a single
active block (for $2\leq k\leq b$),
\item with probability $q(b:1)$ an active block is chosen  uniformly  at random from $b$ blocks and turned into a frozen block.
\end{itemize}
In the case $b=1$ only the second option is possible, that is $q(1:1)=1$, and when 
all blocks of  $\pi^*_n$ are in frozen condition, the operation is defined to be the identity.
For the $\Lambda$-coalescent freezing at positive rate $\rho$, we know that
\begin{itemize}
\item (i) the decrement  matrix  $q$ is of the special form \re{M}, and 
\item (ii)
 the continuous time processes $\Pi^*_n(t)$ are Markovian and consistent as $n$ varies, meaning that $\Pi^*_m(t)$ for
$m <n$ coincides with $\Pi^*_n(t)|_m$, the restriction of $\Pi^*_n(t)$ to $[m]$.
\end{itemize}
Note that  ${\rm FM}_n$ always reduces the number  of active blocks,
in particular it transforms a partition of $[n]$ with $b>1$ active blocks into some other partition of $[n]$
with $b-1$ active blocks with probability $q(b:1)+q(b:2)$.

\par 
To view \Mohle's recursion \re{rec1} in greater generality, we 
consider this freeze-and-merge operation  ${\rm FM}_n$ for $n$ some fixed positive integer, and $q_n$
a finite decrement matrix.
Let $(\widehat{\Pi}_n^*(k), k=0,1,2,\ldots)$ 
be the Markov chain obtained by iterating 
${\rm FM}_n$ starting from $\widehat{\Pi}_n^*(0)=\Sigma_n^*$.  
Since $\FM_n$   is defined in terms of $^*$-shapes,
each $\widehat{\Pi}_n^*(k)$ is a partially frozen exchangeable partition of $[n]$.
The ${\rm FM}_n$-chain is strictly transient, 
in the sense that it never passes through the same state
until it reaches a  
partially frozen partition $E_n^*$, all of whose blocks are frozen.
Let $E_n$ be the induced partition of $[n]$,
 which we call the {\it final partition} and regard $E_n$
as the outcome 
of random transformation
of exchangeable partitions $\Sigma_n\mapsto\Sigma_n^*\mapsto E_n^*\mapsto E_n$.

\par Observe that for $m=1,\ldots, n$ the first $m$ rows of  the decrement matrix $q_n$  comprise a
 decrement matrix $q_m$ which itself defines a freeze-and-merge operation  ${\rm FM}_m$ on partially frozen partitions of $[m]$.
Hence for given $q_n$ we can also define 
a final partition $E_m$ of the ${\rm FM}_m$-chain.
Note that
${\rm FM}_n$ is essentially an  operation on the {\it set} 
of active  blocks, regardless of their contents, sizes, and the 
configuration of frozen blocks. 
\begin{lemma}\label{whyrec1} 
Given an arbitrary decrement matrix $q_n$,
let $p$ be the function on  $\cup_{m=1}^n{\cal C}_m$ whose restriction to 
${\cal C}_m$ is the EPPF of $E_m$, the final partition generated by the
${\rm FM}_m$ chain, for $1\leq m\leq n$.
Then $p$ satisfies \Mohle's recursion {\rm (\ref{rec1})} for  each composition $(n_1,n_2,\ldots,n_\ell)\in {\cal C}_n$.
\end{lemma}
\proof A particular realization of $E_n$ with  $\mathtt{shape}(E_n)=(n_1,\ldots,n_\ell)$ occurs when 
either 
\begin{itemize}
\item (a)
some block $\{j\}$  of $E_n$  appears as a frozen singleton in  ${\rm FM}_n(\Sigma_n^*)$
and all other singletons $\{i\}\neq \{j\}$ evolve to form a partition with shape $(\ldots,\widehat{1},\ldots)$;
or 
\item (b) the first iteration of ${\rm FM}_n$ merges some singletons $\{j_1\},\ldots,\{j_k\}$ ($k>1$) in a single active block
which enters completely one of the blocks of $E_n$.
\end{itemize}
By the definition of $p$ and the last remark before the lemma, the probability of the event (a) is 
$${1\over n}\cdot q(n:1)p(\ldots,\widehat{1},\ldots),$$
because after $\{j\}$ gets frozen the operation ${\rm FM}_n$ is reduced to ${\rm FM}_{n-1}$ acting on partially frozen partitions of 
$[n]\setminus \{j\}$. 
Similarly, the probability of (b) is 
$${1\over{n\choose k}}\cdot q(n:k)p(\ldots,n-k+1,\ldots),$$
because after  creation of the active block $\{j_1,\ldots,j_k\}$ the iterates of ${\rm FM}_n$ can be 
identified with that of
${\rm FM}_{n-k+1}$ acting on partially frozen partitions of $[n]\setminus \{j_2,\ldots,j_k\}$. 
Summation over all possible choices yields (\ref{rec1}).
\endpf

\par In the general setting of Lemma \ref{whyrec1},
the sequence of exchangeable final partitions $(E_m)_{m=1}^n$ need not be consistent with 
respect to restrictions. We turn next to the constraints on $q$ imposed by 
the following stronger consistency condition:

\begin{definition}\label{dd5}
{\rm 
For a decrement matrix $q_{n}$ and $1 \le m < n$, 
call the transition operators ${\rm FM}_n$ and ${\rm FM}_m$  derived from $q_n$
{\it consistent} if whenever $\widehat{\Pi}_n^*$ is a Markov chain governed by ${\rm FM}_n$, the jump process
derived from the restriction of $\widehat{\Pi}_n^*$ to $[m]$ is a Markov chain governed by ${\rm FM}_m$.
Call the decrement matrix $q_{n}$ {\em consistent} if this condition holds for every $1 \le m < n$.
}
\end{definition}

As the leading example, it is clear from consistency of the continuous time chains $(\Pi_n^*(t),~ t \ge 0 )$
which represent a $(\Lambda,\rho)$-coalescent,
that for every $n$ the
corresponding decrement matrix $q_{n}$ is consistent.
The following lemma collects some general facts about consistency. The proofs are elementary and left to the reader.
Let $\FM_n(\pi_n^*)$ denote the random partition obtained by action of $\FM_n$ on an initial
partially frozen partition $\pi_n^*$ of $[n]$,

\begin{lemma}\label{before5} Given a particular  decrement matrix $q_n$:

{\em (i)} For fixed $1 \le m < n$ the transition operators ${\rm FM}_m$ and ${\rm FM}_n$ are consistent if and only if for
each  partially frozen partition $\pi_n^*$ of $[n]$, there is the equality in distribution
$$
\FM_m(\pi_n^*|_{m})\ed\FM_n(\pi_n^*)||_{m}
$$
where on the left side $\pi_n^*|_{m}$ is the restriction of $\pi_n^*$ to $[m]$,
and on the right side the notation $||_{m}$ means the restriction to $[m]$ {\it conditional} 
on the event $\FM_n(\pi_n^*|_{m})\neq \pi_n^*|_{m}$ that $\FM_n$  freezes or merges at least one of the blocks of $\pi_n^*$ 
containing some element of $[m]$.

{\em (ii) } If ${\rm FM}_{m-1}$ and ${\rm FM}_m$ are consistent for every $1 <m \le n$, then 
so are ${\rm FM}_{m}$ and ${\rm FM}_n$ for every $1 <m \le n$; that is, $q_n$ is consistent.

\end{lemma}

\begin{lemma}\label{l5}  A decrement matrix $q_{n}$ is consistent if and only if it satisfies the backward
recursion
\begin{eqnarray}\nonumber
q(b:k)={k+1\over b+1}q(b+1:k+1)+{b+1-k\over b+1}q(b+1:k)~~~~~~~~~~~~~~~~~~~~~~~~~~~~~~\\ \label{rec2}
+{1\over b+1}q(b+1:1)q(b:k)+{2\over b+1}q(b+1:2)q(b:k)~~~~~(2\leq k\leq b<n),\\ \label{rec2b}
q(b:1)={b\over b+1}q(b+1:1)+
{1\over b+1}q(b+1:1)q(b:1)+{2\over b+1}q(b+1:2)q(b:1)~~~(1\leq b<n).\\ \nonumber
\end{eqnarray}
Consequently, each probability distribution $q(n:\cdot)$ on $[n]$ determines a
unique consistent decrement matrix $q_n$ with this $n$th row.
\end{lemma}
\proof
Consider ${\rm FM}_{n}$ and ${\rm FM}_{n-1}$ applied to
$\Sigma_n^*$ and $\Sigma_{n-1}^*$, that is the partitions into singletons, all in the active condition.
For $k\le n-1$, ${\rm FM}_{n-1}$ operates 
by coalescing $\{1,\ldots,k\}$ into an active block with probability 
\eq
\label{30}
\frac{q(n-1:k)}{{n-1 \choose k}}.
\en
As for the jumping process of (${\rm FM}_{n}$ restricted to $[n-1]$),
the probability of a coalescence of $\{1,\ldots,k\}$ into an active block is the sum of the following four parts, depending on the 
development of the ${\rm FM}_{n}$ chain. Let $T_1$ be the time of the first change in the restriction of the
${\rm FM}_{n}$ chain to $[n-1]$.
To obtain the required coalescence, either $T_1 = 1$ and the state after a single step of ${\rm FM}_{n}$ 
comes from $\Sigma_n^*$ by coalescing $\{1,\ldots,k,n\}$ or $\{1,\ldots,k\}$, these occurring with probability
\eq
\label{31}
\frac{q(n:k+1)}{{n \choose k+1}}+\frac{q(n:k)}{{n \choose k}}\, ;
\en
or $T_1 = 2$ with ${\rm FM}_{n}$ acting on
$\Sigma_n^*$ by first freezing $\{n\}$ then coalescing $\{1,2,\ldots,k\}$, or first coalescing $\{n\}$ with one of other $n-1$ 
singletons, leaving $1,2, \ldots k$ in $k$ distinct blocks,  then coalescing these $k$ blocks at the next step; these ways occur with probability
\eq
\label{32}
\frac{q(n:1)}{n}\cdot\frac{q(n-1:k)}{{n-1 \choose k}} + \frac{(n-1)q(n:2)}{{n\choose 2}}\cdot \frac{q(n-1:k)}{{n-1 \choose k}} .
\en
Equate \re{30} with the sum of \re{31} and \re{32} to get (\ref{rec2}) for $b=n-1$.
In much the same way, ${\rm FM}_{n-1}$ may act on
$\Sigma_{n-1}^*$ by freezing $\{1\}$ with probability 
\eq
\label{33}
\frac{q(n-1:1)}{n-1}\,.
\en
While for the  jumping process of (${\rm FM}_{n}$ restricted to $[n-1]$),
to get the required form, either $T_1 = 1$ and ${\rm FM}_{n}$
acts on $\Sigma_n^*$ by freezing $\{1\}$ with probability
\eq
\label{34}
\frac{q(n:1)}{n}\, ;
\en
or $T_1 = 2$ and the result is obtained from $\Sigma_n^*$ by first freezing $\{n\}$ 
then freezing $\{1\}$, or first coalescing $\{n\}$ with one of other $n-1$ singletons then freezing the block containing $1$, 
these ways occurring with probability
\eq
\label{35}
\frac{q(n:1)}{n}\cdot\frac{q(n-1:1)}{n-1}+\frac{(n-1)q(n:2)}{{n \choose 2}}\cdot\frac{q(n-1:1)}{n-1} .
\en
Equate \re{33} with the sum of \re{34} and \re{35} to get (\ref{rec2}) for $b=n-1$.
Combine them to get (\ref{rec2b}) for $b=n-1$. The recursions for $b<n$ follow by replacing $n$ by $b+1$. 

Conversely, granted the recursions (\ref{rec2}) and (\ref{rec2b}), in order to prove consistency
it is enough to check the case $m=n-1$, and this is done by application of Lemma \ref{before5}.
\endpf

\vskip0.5cm

\begin{lemma}\label{l4}
For $1 \le m \le n$ let $E_m$ be the final partition of the $\FM_n$-chain starting in state $\Sigma_m^*$.
If the
decrement matrix $q_n$ is consistent then the finite sequence of exchangeable 
random set partitions $(E_m)_{m=1}^n$ is consistent in the sense that
$$E_m\ed E_n|_m\,.$$
The finite EPPF $p$ of $(E_m)_{m=1}^n$ then satisfies \Mohle's recursion
{\rm (\ref{rec1})} for all 
compositions of $m\le n$ in the left hand side.
\end{lemma}
\proof
The consistency in distribution is clear.
To show (\ref{rec1})
it is enough to look at the case with compositions of $n$ on the left hand side, for which  
Lemma \ref{whyrec1} applies.
\endpf

\par Here is our principal result regarding finite partitions satisfying (\ref{rec1}):

\begin{theorem}\label{thatsit} For a positive integer $n>1$ and arbitrary probability distribution $q(n:\cdot)$ on $[n]$
\begin{itemize}
\item[{\rm(i)}]
there exists a unique finite EPPF $p$ for a consistent sequence of random set partitions $(\Pi_m)_{m=1}^n$ which satisfies \Mohle's recursion {\rm (\ref{rec1})}
for all compositions of $n$ on left hand side, 
\item[{\rm(ii)}]
this finite EPPF $p$ satisfies \Mohle's recursion {\rm(\ref{rec1})} for all compositions of positive integers $m<n$ on the left hand side with coefficients 
$q(m:\cdot)$ derived from $q(n:\cdot)$ by the recursion {\rm (\ref{rec2}), (\ref{rec2b})},
\item[{\rm(iii)}]
for each $1 \le m \le n$ the distribution of $\Pi_m$ determined by the 
restriction of this EPPF $p$ to compositions  of $m$ is that of
the final partition of the $\FM_m$ Markov chain with 
 decrement 
matrix $q_m$ defined  by {\rm (ii)}, 
starting from state $\Sigma_m^*$.

\end{itemize}
\end{theorem}

\proof We apply Lemma \ref{l4}.
Given arbitrary probability distribution $q(n:\cdot)$ on $[n]$, we can define all $q(m:\cdot)$, $1\le m<n$, by the backward recursion (\ref{rec2}),
(\ref{rec2b}).
Then we use the decrement matrix $q_n$ with these rows
to build a sequence of Markov chains: for each $m$, 
the chain $(\Pi_m(k),k=0,1,2,\ldots)$ starts from $\Sigma_m^*$ and evolves according to $\FM_m$. 
The  sequence of induced final partitions $(E_m)_{m=1}^n$ of these chains has
 EPPF $p$ which satisfies recursion (\ref{rec1}). Hence the existence part of (i) follows.
We postpone the proof of uniqueness in part (i) to the next section.
The assertions (ii) and (iii) follow
directly from this  construction.
\endpf

\section{The sample-and-add operation} \label{5}

Given a probability distribution $q(n:\cdot)$ on $[n]$, we now interpret
\Mohle's recursion {\rm(\ref{rec1})} as the system of equations
for the invariant probability measure of a particular Markov transition
mechanism on partitions of $[n]$, and show that this invariant probability
distribution is unique.
This will complete the proof of Theorem \ref{thatsit}.

\par Consider the following {\em sample-and-add} random operation on ${\cal P}_{[n]}$,
denoted ${\rm SA}_n$.
We regard a generic random partition $\Pi_n\vdash [n]$ as a random allocation 
of balls labeled $1,\ldots,n$ to some set of nonempty boxes,
which the operation ${\rm SA}_n$ transforms into some other 
random allocation $\Pi_n'$.
 Fix $q(n:\cdot)$, 
a probability distribution on $[n]$
and let
$K_n$ be a random variable with this distribution $q(n:\cdot)$. 
Given $K_n = k$ and $\Pi_n = \pi_n$,
\begin{itemize}
\item
if $k=1$, first delete a single ball picked uniformly at random from the balls
allocated according to $\pi_n$, to make an intermediate 
partition of some set of $n-1$ balls, then add to this intermediate 
partition a single box containing the deleted ball.
\item if $k=2,\ldots,n$, delete a sequence of $k-1$ of the $n$ balls from
$\pi_n$ by uniform random sampling without replacement, to obtain
an intermediate partition of some set of $n-k+1$ balls, then mark a ball 
picked uniformly from these $n-k+1$ balls, and add the $k-1$ sampled balls 
into the box containing the marked ball. 
\end{itemize}
In either case delete empty boxes in case any appear after the sampling
step.  The resulting partition of $[n]$ is $\Pi_n'$. 
For each $q(n:\cdot)$, this defines a Markovian transition operator
${\rm SA}_n$ on partitions of $[n]$.

\begin{lemma}\label{stat} 
Let $\Pi_n$ be an exchangeable random partition of $[n]$ with finite
EPPF $p$ defined as a function of compositions of $m$ for $1 \le m \le n$.
Let $\Pi_n'$ be derived from $\Pi_n$ by the ${\rm SA}_n$ operation determined
by some arbitrary probability distribution $q(n:\cdot)$ on $[n]$.
Then $\Pi_n'$ is an exchangeable random partition of $[n]$ whose EPPF
$p'$ is determined on compositions of $[n]$ by the formula
\eq\label{rec11}
p'(n_1,n_2,\ldots,n_\ell)=
 \frac{q(n:1)}{n}\sum_{j:n_j=1} p(\ldots,\widehat{n_j},\ldots)+ 
\sum_{k=2}^n q(n:k)\sum_{j: n_j\ge k} \frac{{n_j\choose
k}}{{n\choose k}}p(\ldots, n_j-k+1,\ldots).
\en
\end{lemma}

\noindent
{\em Note.} The right side of \re{rec11} is identical to the right side of
\Mohle's recursion (\ref{rec1}).\\
\proof
Let $K_n$ with distribution $q(n:\cdot)$ be the number of balls
deleted in the ${\rm SA}_n$ operation. For each partition $\pi_n'$ 
of $[n]$ we can compute
\eq\label{rec111}
\prob ( \Pi_n' = \pi_n') = \sum_{k=1}^n q(n: k)  \,
\prob ( \Pi_n' = \pi_n' \giv K_n = k).
\en
Assuming that $\pi_n'$ has boxes of sizes $n_1, \ldots, n_\ell$,
and that the ${\rm SA}_n$ operation acts on an exchangeable $\Pi_n$ with 
EPPF $p$, we deduce \re{rec11} from \re{rec111} and
\eq
\label{rec1one}
\prob ( \Pi_n' = \pi_n' \giv K_n = 1) = \frac{1}{n}\sum_{j:n_j=1} p(\ldots,\widehat{n_j},\ldots), 
\en
\eq
\label{rec1111}
\prob ( \Pi_n' = \pi_n' \giv K_n = k) =  \sum_{j: n_j\ge k}  
\frac{{n_j\choose k}}{{n\choose k}}p(\ldots, n_j-k+1,\ldots) , \qquad k \ge 2.
\en
Consider \re{rec1111} first.
For the event $(\Pi_n' = \pi_n')$ to occur there must be some $j$ 
with $n_j \ge k$.
For each such $j$, corresponding to a box of $\pi_n'$ with at least
$k$ balls, the result $(\Pi_n' = \pi_n')$ might be obtained by addition
of $k-1$ balls to that box.
The sequence of labels of these balls, in order of their choice, can
be any one of
$
n_j (n_j-1) \cdots (n_j - k + 2)
$
sequences, and the final ball chosen to mark the box can be any one of
$n_j - k + 1$ balls, making $k! { n_j \choose k }$ choices out of a 
total of $k! { n \choose k }$ possible choices.
Given one of these $k! { n_j \choose k }$ choices of $k$ balls, let $M_{k-1}$
be the set of labels of the $k-1$ balls that are moved.
Then the event $(\Pi_n' = \pi_n')$ occurs if and only if 
the restriction of $\Pi_n$ to $[n] - M_{k-1}$ equals the
restriction of $\pi_n'$ to $[n] - M_{k-1}$, which is a particular partition
of $n-k+1$ labeled balls into boxes of $\bar{n}_1,  \ldots, \bar{n}_\ell$ 
balls, where $\bar{n}_i = n_i 1(i \ne j ) + (n_j - k + 1) 1(i = j)$.
The conditional probability of $(\Pi_n' = \pi_n')$,
given $K_n = k$ and which of the $k! { n_j \choose k }$ possible
choices of $k$ balls is made, is therefore  $p(\ldots, n_j -k + 1, \ldots)$,
by the assumed exchangeability of $\Pi_n$, and the definition of the
EPPF $p$ of $\Pi_n$ on compositions of $m \le n$ by restriction of
$\Pi_n$ to subsets of size $m$.
The evaluation  \re{rec1111} is now apparent, and \re{rec1one} too is
apparent by a similar but easier argument.

\endpf

\begin{proposition}\label{sainv}
For each probability distribution $q(n:\cdot)$ on $[n]$, the corresponding
${\rm SA}_n$ transition operator on partitions of $[n]$ 
has a unique stationary distribution.
A random partition with this stationary distribution is
exchangeable, and its EPPF is the finite unique EPPF $p$ that
satisfies \Mohle's recursion {\rm (\ref{rec1})}, that is
{\rm \re{rec11}} with $p'=p$.
\end{proposition}

\proof
If $q(n:1)=1$ then eventually ${\rm SA}_n$ terminates with singleton 
partition, 
so the stationary distribution is degenerate and concentrated on the singleton partition.
If $q(n:1)=0$ then eventually ${\rm SA}_n$ terminates with one-block 
partition, 
so the stationary distribution is degenerate and concentrated on the one-block
partition.
If $0<q(n:1)<1$ then also $q(n:k)>0$ for some $k>1$; in this 
case 
the stationary  law is again unique because all states communicate:
e.g. the pure-singleton partition $\Sigma_n$ is reachable from everywhere, and it can reach any partition in finitely many steps, as is easily 
verified.
Observe that passing to shapes projects 
the ${\rm SA}_n$ chain with state
space partitions of the set $[n]$ onto another Markov chain whose state space 
is the set of partitions of the integer $n$. It follows easily that
the unique stationary distribution of
${\rm SA}_n$ governs an exchangeable random partition of $[n]$.
The previous lemma shows that its EPPF $p$ solves
\Mohle's recursion.
Finally, if an EPPF $p$ solves \Mohle's recursion, then it provides a
stationary state for the ${\rm SA}_n$ chain. Hence the uniqueness result for
solutions of \Mohle's recursion by an EPPF $p$.
\endpf

\subsection{Special cases}
Following are two special cases of ${\rm SA}_m$ operation:

\vskip 0.5cm

\noindent \textbf{Ewens' partition} appears when $q(n:\cdot)$  may have only two positive  entries
$$q(n:1)={2\rho \over n-1+2\rho}\,\,~~~{\rm and~~}~q(n:2)={n-1\over n-1+2\rho}$$
for each $n\ge 2$. It is easy to realize that the ${\rm SA}_n$ operation in this case is reduced to the following operation
 with $u=2\rho / (n-1+2\rho)$: given a number $0\le u\le 1$ and a partition of $[n]$ as allocation of $n$ labeled balls, first we 
uniformly sample two balls named $A$ and $D$ without replacement from the $n$ balls (so $A = D$ is excluded), then we 
put ball $A$ back to where it was, and finally
 \begin{itemize}
 \item with probability $u$ append a new box containing the single ball $D$,
 \item with probability $1-u$ add the ball $D$ to the box containing ball $A$.
\end{itemize}
In this case, if we consider the FM operator determined by $q$, it is clear that only binary merges happen. 
That the stationary partition $\Pi_n$ follows the Ewens' sampling formula with  parameter $\theta=(n-1)u/(1-u)$ is seen
by the `Chinese restaurant' rule \cite{CSP} for transition from $\Pi_{n-1}$ to $\Pi_n$, or
can be easily concluded from  the formula. The coincidence 
of the stationary distribution of this ${\rm SA}_n$ chain with the law of the induced final 
partition $E_n$ of the associated $\FM_n$ chain confirms in this case the well known fact
that Kingman's coalescent with mutations terminates at Ewens' partition.

\par
The ${\rm SA}_n$-chain resembles
Moran's novel mutation  chain \cite{moran58,  watterson76, young95}. 
Transitions of the latter 
are the  following: 
given a number $0\le u\le 1$ and a partition of $[n]$ as allocation of labeled balls, 
first choose two balls named $A$ and $D$ uniformly and independently from the $n$ balls (so $A=D$ is not excluded), then
follow the rules
\begin{itemize}
\item with probability $u$ append a new box with a single ball $C$,
\item with probability $1-u$  add a ball $C$ to the box that  contains  ball $A$, 
\end{itemize}
then  assign to ball $C$ the same label as that of $D$ and 
finally remove ball $D$. 
It is well known \cite{watterson76} that the stationary law of Moran's  chain 
corresponds to Ewens' partition  with parameter $nu/(1-u)$.
\vskip0.5cm

\noindent \textbf{Hook partitions.} Another extreme case appears when $q(n:\cdot)$ may have  only two positive entries 
$$q(n:1)={n\rho\over 1+n\rho}\,\,~~~{\rm and}~~~q(n:n)={1\over 1+n\rho}.$$
In this case ${\rm SA}_n$ creates some number of singletons and then after some number of steps
puts all balls in a single box. If $0<q(n:1)<1$,
the stationary distribution concentrates on partitions with a hook shape  $(m,1,1,\ldots,1)$.
This partition results from the $\Lambda$-coalescent with freeze when 
$\Lambda=\delta_1$ is a Dirac mass at $1$.

\section{Infinite partitions}
\label{6}

In this section we pass from finite partitions to the projective limit, and
arrive at the 
desired integral representation of infinite decrement matrix $q_\infty$
satisfying recursion (\ref{rec2}), (\ref{rec2b}). 
This 
gives another approach 
to \Mohle's partitions via consistent freeze-and-merge chains, which may be seen as  discrete-time jumping 
processes associated with the $\Lambda$-coalescent with freeze.

\par An infinite sequence of freeze-and-merge operations $\FM:=(\FM_n, n=1,2,\ldots)$ 
which satisfies the condition in Definition \ref{dd5} for all positive integers $1\le m<n<\infty$ is called {\it consistent}. 
By Lemma \ref{l5} such a sequence $\FM$ is determined by
an infinite decrement matrix $q_\infty$ which satisfies the recursions
 (\ref{rec2}), (\ref{rec2b}). 

 \par For each $n=1,2,\ldots$
the Markov chain starting from $\Sigma_n^*$ and driven by $\FM_n$ terminates with an induced final partition $\Pi_n$.
These comprise an infinite partition $\Pi_\infty=(\Pi_n)_{n=1}^\infty$ which we  call the {\it final partition} associated with consistent $\FM$.
In the case $q(2:1)=0$ the final partition is the trivial one-block partition.

\begin{lemma}\label{l6} For every infinite decrement matrix $q_\infty$ with entries satisfying the recursion {\rm(\ref{rec2}), (\ref{rec2b})} 
there exist a non-negative finite measure $\Lambda$ on $[0,1]$ 
and a non-negative real number $\rho$, which satisfy $(\Lambda,\rho)\neq (0,0)$ and are such that  the representation
$q(n:k)=\Phi(n:k)/\Phi(n)$  $(1\leq k\leq n)$  holds 
with $\Phi$ as in \,{\rm (\ref{Phi1})},
{\rm (\ref{Phi2})}, {\rm (\ref{Phi3})}.
The data 
$(\Lambda,\rho)$ are unique up to a positive factor.
\end{lemma}
\proof
Suppose $q$ solves (\ref{rec2}), (\ref{rec2b}) and suppose $q(2:2)<1$.
Let $\Phi(n), n=1,2,\ldots$ satisfy
\eq\label{PP}
{\Phi(n)\over\Phi(n+1)}=1-{1\over n+1}q(n+1:1)-{2\over n+1}q(n+1:2)
\en
for $n\geq 1$; because the right side is strictly positive this recursion has a unique  
solution with some given initial value $\Phi(1)=\rho$, where $\rho>0$.
For $2\leq k\leq n$ set
$$\Phi(n:k):=q(n:k)\Phi(n),$$ 
then 
from (\ref{PP}) and (\ref{rec2})
$$\Phi(n:k)=
{k+1\over n+1}\Phi(n+1:k+1)+{n+1-k\over n+1}\Phi(n+1:k)\qquad~~~~(2\leq k\leq n<\infty).$$
Apart from a shift by $2$, this is the well-known Pascal-triangle recursion appearing in connection 
with de Finetti's theorem and the Hausdorff moment problem, hence
(\ref{Phi2}) holds for some non-negative measure $\Lambda$ on Borel sets of $[0,1]$. 
From (\ref{rec2b}) we find
$$\rho={\Phi(1)q(1:1)\over 1 }=\cdots={\Phi(n)q(n:1)\over n}=\cdots,$$
and from 
$$\sum_{k=1}^n\Phi(n)q(n:k)=\Phi(n)$$
we deduce (\ref{Phi3}) and $q(n:1)=\rho n/\Phi(n)$. Setting by definition $\Phi(n:1):=\rho n$ we are done.
For the special case $q(2:2)=1$, it is easy to observe that $\rho =0$, and we get $\Lambda=\delta_0$ 
by similar analysis.
\endpf

\vskip0.5cm

\par Recording this lemma together with previous results, we have:

\begin{theorem} \label{MAIN}
Let $(\Pi_n)_{n=1}^\infty$ be a nontrivial exchangeable random partition of $\Nat$, 
different from the trivial one-block partition.
The following are equivalent:
\begin{itemize}
\item[{\rm(i)}] The EPPF $p$ satisfies recursion {\rm (\ref{rec1})} with some infinite decrement matrix 
$q_\infty$.
\item[{\rm(ii)}] This matrix is representable as 
$q(n:k)=\Phi(n:k)/\Phi(n)$
with $\Phi$ defined by 
{\rm (\ref{Phi1}), (\ref{Phi2}), (\ref{Phi3})} and
some nontrivial $(\Lambda,\rho)$, which is   unique up to a positive factor.
\item [{\rm(ii)}] This $\Pi_\infty$ is induced by the final partition of some standard  $\Lambda$-coalescent freezing at rate $\rho$.
\item[{\rm (iii)}] This $\Pi_\infty$ is the final partition of some consistent $\FM$ operation.
 \end{itemize}
\end{theorem}

\vskip0.5cm
\par Complementing this result, we have the following uniqueness assertion.

\begin{lemma}\label{inj}
The correspondence $q\mapsto p$ between infinite decrement matrices with $q(2:1)>0$ satisfying consistency {\rm  (\ref{rec2}), (\ref{rec2b})}
 and the EPPF's is bijective.  
\end{lemma}
\proof
 We only need to show that $p$, which by Lemma \ref{l4} must solve (\ref{rec1}), uniquely determines $q$.
For general infinite partitions $q(2:1)=p(1,1)>0$  
implies that $p(1,1,\ldots,1)>0$.
This applied to the singleton shapes together with
$$p(1,\ldots,1)=q(n:1)q(n-1:1)\cdots q(2:1)$$
shows that the $q(n:1)$'s  are uniquely determined by $p$. To show that
$q(n:m)$ for $1\le m\le n-1$ is also determined by $p$,
exploit the formula 
\begin{eqnarray*}
p(m,1,\ldots,1)={q(n:m)\over{n\choose m}}p(1,\ldots,1)+~~~~~~~~~~~~~~~~~~~~~~~~~~~~~~~~~~~~~~~~ \\
\sum_{k=2}^{m-1}
q(n:k)\frac{{m\choose k}}{{n\choose k}}p(m-k+1, 1,\ldots,1)+
q(n:1)\frac{n-m}{n}p(m,\widehat{1},1,\ldots,1),
\end{eqnarray*}
and argue by induction in  $m=2,3,\ldots,n-1$.
\endpf
\vskip0.5cm

Thus if an exchangeable infinite partition can be realized 
as the induced final partition of a consistent {\rm FM}-operation, then this {\rm FM}-operation is unique.
The realization via a $(\Lambda,\rho)$-coalescent process is unique up to a positive multiple of the parameters, 
which corresponds to a linear time-change of the coalescent.
If there is no freeze the uniqueness fails, since any $\Lambda$-coalescent terminates with the trivial one-block partition.

\par We classify next the   cases when some of the entries of $q$ are zeros. It is assumed that the starting partition is $\Sigma_\infty^*$.  
\begin{itemize}

\item[(i)] If $q(n:1)=1$ holds for $n=2$ then the same holds for $n\geq 2$. This is the pure-freeze coalescent with $\Lambda=0$, hence $E_\infty=\Sigma_\infty$.

\item[(ii)] If $q(n:1)=0$ holds for $n=2$ then the same holds for $n\geq 2$. This is a $\Lambda$-coalescent with no freeze, 
hence $E_\infty$ is the one-block partition.

\item[(iii)] If $q(n:1)>0, q(n:2)>0$ and  $q(n:1)+q(n:2)=1$ hold  for $n=3$ then the same relations hold for $n\geq 3$.
This is the case of Kingman's coalescent with freeze, $\Lambda$ is a positive mass at $0$, and $E_\infty$ is Ewens' partition. 
ld

\item[(iv)] if $q(n:1)>0, q(n:n)>0$ and $q(n:1)+q(n:n)=1$ hold 
for $n=3$ then also for $n\geq 3$.
In this case 
$\Lambda$ is a positive mass at $1$, and $E_\infty$ is a  hook partition.

\end{itemize}
\noindent
The `generic' case is characterised by 
$q(3:1)>0, q(3:2)>0, q(3:3)>0$, in which case $0<q(n:m)<1$ for all $1\leq m\leq n<\infty$. 

\section{Positivity} \label{pos}

This section provides a construction of  decrement matrices
$q_\infty$ satisfying the consistency condition (\ref{rec2}), (\ref{rec2b}),
from a single sequence of real numbers satisfying a positivity condition.
For  $(c(n), n=0,1,2,\ldots)$ a  sequence of real numbers, the backward difference operator $\nabla$ is defined as
$$
\bigtriangledown c(n):=c(n)-c(n+1),
$$
and for any $j=0,1,2,\ldots$ its iterates act as 
$$
\bigtriangledown^j c(n)=\sum_{i=0}^j(-1)^i {j\choose i}c(n+i).
$$

Now let $(\Phi(n), n=1,2,\ldots)$ be a sequence of real numbers and $\rho$ be a positive real number. Define for each $n$
\eq\label{phin1}
\Phi(n:1) := \rho n,
\en
and
\eq
\overline{\Phi}(n):=\Phi(n)-\rho n \,.
\en
Define 

\eq\label{psi}
\Psi(n):=\frac{\bigtriangledown\overline{\Phi}(n)}{n}
\en
and let 
\eq\label{phinm}
\Phi(n:m):=-{n \choose m}\bigtriangledown^{m-2}\Psi(n-m+1),\ \ \ 2\le m\le n.
\en

With these definitions, it can be verified that for each $n$ 
\eq
\Phi(n)=\Phi(n:1)+\Phi(n:2)+\cdots+\Phi(n:n).
\en
Hence if all $\Phi(n)$ are positive and all $\Phi(n:m)$ are non-negative, the matrix with entries
\eq\label{decrementfromphi}
q(n:m):=\frac{\Phi(n:m)}{\Phi(n)},\ \ \ \ 1\le m\le n 
\en
is a well defined infinite decrement matrix. More than that, we have the following observation:
\begin{lemma}
Suppose that a sequence of positive real numbers $\rho$, $\Phi(n)$, $n=1,2,\ldots$ is such that each entry $\Phi(n:1)$, $\Phi(n:m)$ 
in {\rm (\ref{phin1}), (\ref{phinm})} is non-negative. Then the matrix {\rm (\ref{decrementfromphi})} 
satisfies the recursion {\rm (\ref{rec2}), (\ref{rec2b})}.
\end{lemma}
\proof
The definition (\ref{phinm}) of $\Phi(n:m)$ implies the recursion
\eq\label{triang}
\Phi(n:m)=\frac{m+1}{n+1}\Phi(n+1:m+1)+\frac{n-m+1}{n+1}\Phi(n+1:m),\ \ \ 2\le m\le n.
\en
Using this relation, the first recursion (\ref{rec2}) can be reduced to 
$$
2\Phi(n+1:2)=(n+1)(\Phi(n+1)-\Phi(n))-\Phi(n+1:1)
$$
which follows from definition of $\Phi(n+1:2)$ and $\Phi(n+1:1)$.
The second recursion is actually the definition of $\Phi(n+1:2)$ after we plug in all the $\Phi(n:1)$, $\Phi(n+1:1)$ terms.
\endpf
\vskip0.3cm

The above lemma shows that given a sequence of positive real numbers with some additional positivity property, we can 
recover \Mohle's partition structure by first defining a consistent decrement matrix, then using the recursion (\ref{rec1}). 
By Lemma \ref{l6} we know that every decrement matrix satisfying consistency condition (\ref{rec2}) (\ref{rec2b}) has an integral
 representation which is unique up to a positive factor, so it is clear that we also have integral representation for the sequence 
of $\Phi(n)$ given here:
\begin{proposition}
A sequence of positive real numbers $\rho$, $\Phi(n)$, $n=1,2,\ldots$ is such that each entry $\Phi(n:1)$, 
$\Phi(n:m)$ as in {\rm (\ref{phin1}) , (\ref{phinm})} is non-negative if and only if these numbers admit the integral representation 
{\rm \rem{Phi1},\rem{Phi2},\rem{Phi3}} for some non-negative finite 
measure $\Lambda$ on $[0,1]$, which is then unique.
\end{proposition}

\section{Freezing times }
\label{7}

In this section $(\Pi^*(t), t\geq 0)$ is a standard $(\Lambda,\rho)$-coalescent, with
$(\Pi(t), t\geq 0)$ induced ordinary partitions, and $E_\infty$ final partition.
We assume that both $\Lambda$ and $\rho$ are nonzero.
The process $(\Pi^0(t), t\geq 0)$ will denote the standard $\Lambda$-coalescent.
We presume that all $(\Lambda,\rho)$-coalescents are defined consistently as $\rho$ varies,
so that the $\Pi(t)$'s and $E_\infty$ get finer as  the freezing rate $\rho$ increases, in particular  each partition $\Pi(t)$ is finer than 
$\Pi^0(t)$, for each $t\geq 0$ and $\rho>0$.

\subsection{Age ordering}

Assigning each individual $j\in \Nat$ the {\it freezing time} $\tau_j$, when the active block containing $j$ 
gets frozen, the final partition $E_\infty$ is defined by sending $i,j$ to the same block if and only if
$\tau_i=\tau_j$. The correspondence $j\mapsto\tau_j$ induces a total order on the set of blocks of $E_\infty$: we say that
the block containing $j$ is {\it older} than the block containing $i$ if $\tau_i<\tau_j$.
With this {\it age ordering}, $E_\infty$ is an {\it ordered exchangeable partition} of $\Nat$, as
studied in \cite{donnellyjoyce91,donnellytavare86,RCS}. 

\par We preserve the notation $E_\infty=(E_n)$ to denote the partition with this additional feature of total order
on the set of the blocks. The law of ordered partition $E_\infty$ is determined by an exchangeable 
composition probability function (ECPF) $c(n_1,\ldots,n_\ell)$ on compositions of $n$.
The ECPF $c$ must satisfy an addition rule similar to (\ref{sps}) but, unlike $p$, need not be symmetric.
The EPPF $p$ of unordered partition is recovered from $c$ by symmetrization. See \cite{RCS} for details.

\par With each $j$ we associate a random open interval $]a_j, b_j[\,$, where
\eq\label{ab}
a_j=\lim_{n\to\infty} \#\{i\leq n: \tau_i<\tau_j\}/n\,,~~~b_j-a_j=\lim_{n\to\infty} \#\{i\leq n: \tau_i=\tau_j\}/n\,,
\en
and the existence of the frequencies is guaranteed by de Finetti's theorem. Thus $a_j$ is the total frequency of blocks 
preceding the block containing $j$, and $b_j-a_j$ is the frequency of the block containing $j$.
The random open set $ U=\cup_j \,]a_j,b_j[$ is the {\it paintbox} representing $E_\infty$.
The partition $E_\infty$ can be uniquely recovered from $U$ by a simple sampling scheme \cite{kingman78a,kingman82a,RCS}.

\par For instance, when $\Lambda=\delta_0$, the complement closed set is $U^c=\{1,\, Y_1,\, Y_1Y_2,\, \ldots,\, 0\}$ for $Y_k$'s 
independent random variables whose distribution is beta$(2\rho,1)$. This case has been thoroughly studied 
\cite{donnellyjoyce91,donnellytavare86},
and it is well 
known that the arrangement of the block sizes in the age order is inverse to the arrangement in size-biased order.
In the case $\Lambda=\delta_1$, the set $U$ has only one interval $]Y,1[$, where $Y$ has a beta distribution.

\subsection{Properties of the final partition}

Some properties of $U$ for  a $(\Lambda,\rho)$-coalescent with $\rho>0$ follow from known results about the $\Lambda$-coalescents
\cite{lambda}. We shall discuss only the case $\Lambda\{1\}=0$, since the case $\Lambda\{1\}>0$ 
only differs by an independent exponential killing and its properties easily follow from that in the  case $\Lambda\{1\}=0$.
Let
$$\mu_r:=\int_0^1 x^r\Lambda(d x).$$
Denote $\rm Leb$ the Lebesgue measure on $[0,1]$.
In the event ${\rm Leb}(U)<1$ the ordered partition $E_\infty$ with paintbox $U$ has a positive total frequency of singletons
blocks, and in the event ${\rm Leb}(U)=0$ there are no singleton blocks at all.

\begin{proposition}\label{infsingl} If $\mu_{-1}<\infty$ then with probability one 
\begin{itemize}
\item[\rm (i)] $\Pi^0(t)$ has singletons, for each $t>0$,
\item[\rm (ii)] $\Pi^*(t)$ has active singletons, for each $t>0$,
\item[\rm (iii)] $\Pi^*(t)$ has frozen singletons, for each $t>0$,
\item[\rm (iv)]
$E_\infty$ has singleton blocks.
\end{itemize}
If $\mu_{-1}=\infty$ then the opposites of {\rm (i)-(iv)} hold with probability one.
\end{proposition}
\proof
By \cite[Lemma 25]{lambda}, if $\mu_{-1}<\infty$ then $\Pi^0(t)$ has singletons almost surely , and if  $\mu_{-1}=\infty$
 the partition has no singletons almost surely.
Now,  if 
 $\Pi^0(t)$ has singletons each of them is active  with probability $0<e^{-\rho t}<1$, independently of the others, thus
the partially frozen partition $\Pi^*(t)$ has singletons in both conditions, and the frozen ones   are also 
singleton blocks of $E_\infty$.
Conversely, if with positive probability $E_\infty$ has singletons then for some $t$ with positive probability $\Pi^*(t)$ has frozen singletons,
then,  perhaps for some other $t$, with positive probability   $\Pi^*(t)$ has active singletons,
but in this event the partition $\Pi^0(t)$ has singletons, hence $\mu_{-1}=\infty$ cannot hold.
\endpf

\par By \cite[Proposition 23]{lambda}   
the $\Lambda$-coalescent either {\it comes down from  infinity}
(the number of blocks in $\Pi^0(t)$, is finite almost surely
for every $t>0$)
or  {\it stays infinite}  (the number of blocks is finite).

\begin{proposition}
If the 
 $\Lambda$-coalescent stays infinite, then the $(\Lambda,\rho)$-coalescent has infinitely many active blocks at any time, therefore
\begin{itemize}
\item[\rm (i)] the set of freezing times $\{\tau_j\}$ is dense in ${\mathbb R}_+$,
\item[\rm (ii)] the closed set $U^c$ has  empty interior and no isolated points.
\end{itemize}
If  the 
 $\Lambda$-coalescent comes down from  infinity, then the $(\Lambda,\rho)$-coalescent satisfies
\begin{itemize}
\item[\rm (i$'$)] the set of freezing times $\{\tau_j\}$ is bounded and only accumulates near $0$,
\item[\rm (ii$'$)] the closed set $U^c$ only accumulates near $0$.
\end{itemize}
\end{proposition}
\proof
Let $J_k$ be the minimal element in some block $A_k$ of $\Pi^0(t)$. Then $J_k$ is also the minimal element in some block 
$B_k\subset A_k$ of $\Pi^*(t)$. Since the block containing  $J_k$ changes the condition from active to frozen  independently of
the $\Lambda$-coalescent,
with positive probability $1-e^{-\rho t}$ the block $B_k$ is active.
For $k=1,2,\ldots$ these events are independent, hence $\Pi^*(t)$ has infinitely many active blocks.
But the same is true for $t+\epsilon$, hence arguing as in Proposition \ref{infsingl} we see that infinitely
many of the active $B_k$'s get frozen before $t+\epsilon$, whence (i). Moreover,
infinitely many of the active $B_k$'s are nonsingleton, hence, by
the law of large numbers for exchangeable trials, have positive frequency.
The assertion (ii) follows now from this remark, (i) and (\ref{ab}). 
\endpf

\section{Comparision with regenerative partitions} \label{8}

This section is devoted to parallels and differences between 
\Mohle's partitions and regenerative partitions
\cite{RCS, RPS}. A novel feature discussed here is a realization of
regenerative partitions by a simple continuous-time coalescent process.

\subsection{Continuous time realization and EPPF}

Consider a ${\cal P}^*_\infty$-valued Markovian process $(\Pi_\infty^*(t),\,t\geq 0)$ which starts 
with $\Sigma_\infty^*(t)$ and evolves by the following
rules.
Any number of 
{\it active singleton} blocks  can merge to form a {\it single frozen} block,
which suspends further evolution immediately. In particular,
an active singleton block can turn into frozen singleton block, an event interpreted as unary merge.
If $\Pi_n(t)$ has $b$ active blocks, each $k$-tuple is merging at the same rate, so that the total
rate for a $k$-merge is $\Phi(b:k)$, for $1\leq k\leq b<\infty$, and
$\Phi(1:1)>0$.

\par Eventually there are only frozen blocks whose configuration determines a final partition   $E_\infty$.
Setting $\Phi(b):=\Phi(b:1)+\ldots+\Phi(b:b)$ and $q(n:k):=\Phi(n:k)/\Phi(n)$,
the EPPF of $E_\infty$ satisfies
\eq\label{EPFRegen}
p(n_1,n_2,\ldots,n_\ell)=
\sum_{j=1}^\ell \frac{1}{ {n \choose n_j } }\, q(n:n_j)\, p(\ldots,\widehat{n_j},\ldots)
\en
for any composition $(n_1,n_2,\ldots,n_\ell)$ of $n$, which is 
a recursion analogous to (\ref{rec1}).
This allows an explicit formula
\eq
\label{recregen1}
p(n_1,n_2,\ldots,n_\ell)=
\sum_{\sigma} \,{ q(N_{\sigma(1)}:n_{\sigma(1)})\cdots q(N_{\sigma(\ell)}:n_{\sigma(\ell)})
\over {n \choose n_1,\ldots,n_\ell}}
\, .
\en
where the sum is over all permutations $\sigma:[\ell]\to[\ell]$, and $N_{\sigma(j)}=n_{\sigma(j)}+\ldots+n_{\sigma(\ell)}$.

\subsection{Subordinator}

 Exchangeability implies the existence of a nonnegative 
finite measure on $[0,1]$  such that
\eq\label{PhiRegen}
\Phi(b:k)={n\choose k}\int_0^1 x^{k-1}(1-x)^{b-k}\Lambda({
\rm d}x),
\en
a representation
to be compared with (\ref{Phi2}).
The cumulative rate for some transition  when  $\Pi_n(t)$ has $b$ active blocks equals 
$$
\Phi(b):=\Phi(b:1)+\ldots+\Phi(b:b)=\int_0^1 {1-(1-x)^b\over x}\,\Lambda({\rm d}x).
$$

\par The last formula  is an integral representation of a Bernstein function, hence 
the measure $\Lambda({\rm d}x)/x$ can be associated with some subordinator  \cite{gnedinpitman06}.
Explicitly, by de Finetti's theorem there exists 
 the limit proportion $S_t$ of integers in $[n]$ that comprise the active blocks of  $\Pi_n^*(t)$, as $n\to\infty$.
The process $(-\log(1-S_t),\,t\geq 0)$ is a subordinator
with $S_0=0$ and distribution determined by 
$$
\ex [(1-S_t)^\lambda]= e^{-t\Phi(\lambda)},~~~~~~t\geq 0,\,\lambda\geq 0,
$$ 
which is a 
version of the \Lev-Khintchine formula in the form of the Mellin transform.
The subordinator has a drift if $\Lambda$ has an atom at $0$.

\par Putting the blocks of $E_\infty$ in increasing order of their  freezing times yields an ordered exchangeable partition
with ECPF 
$$p(n_1,\ldots,n_\ell)=\prod_{j=1}^\ell {q(N_j:n_j)\over {N_j\choose n_j}}\,,$$
where $N_j:= n_j + \cdots + n_\ell$.
The closed range of the process $(S_t)$ is the complement $U^c$ to the paintbox $U$ of the ordered partition $E_\infty$.

\subsection{Related Markov chains }

\subsubsection{Transient}

For regenerative partitions the analogue of  $\FM_n$ introduced in Section \ref{4} is the following.
Let $q_\infty=\{q(b:k), 1\le k\le b<\infty\}$ be a decrement matrix.
If there are $b$ active blocks in a partially frozen partition of $[n]$, then with probability $q(b:k)$ any $k$ of $b$ active blocks
 are chosen uniformly at random and merged into a single frozen block. 
\par Consistency translates as the recursion 
\eq\label{recq}
q(b:k)=\frac{k+1}{b+1}q(b+1:k+1)+\frac{b+1-k}{b+1}q(b+1:k)+\frac{1}{b+1}q(b+1:1)q(b:k)
\en
with $q(1:1)=1$, which leads to 
$$
q(b:k)=\Phi(b:k)/ \Phi(b)~~~~~(1\le k\le b<\infty),
$$
where $\Phi$ has the above integral representation (\ref{PhiRegen}) 
with some measure $\Lambda$ unique up to a positive multiple.

\subsubsection{Recurrent}

The analogue of operation 
${\rm SA}_n$ introduced in Section \ref{5}, acting on ordinary partitions of $[n]$,  is the following \cite{RPS}. 
 Given a decrement matrix $q$, 
let $K_n$ follow $q(n:\cdot)$. Choose a value $k$ for $K_n$, then starting from some partition $\pi_n$ of $[n]$
sample $k$  balls from $\pi_n$ uniformly without replacement, and then append a new box with these $k$ balls to the remaining partition
of $n-k$ balls.
According to an ordered version of the algorithm, acting on ordered partitions,
the balls are sampled from a totally ordered series of boxes, and the newly created box is always 
arranged as the first box in the series.

\par In contrast to the ${\rm SA}_n$ operation, these Markov chains on partitions 
of $[n]$ are consistent under restrictions as $n$ varies. 
To see that {\it the operations ${\rm SA}_n$ are not consistent as $n$ varies} (exluding the hook case $q(n:1)+q(n:n)\equiv 1$)  fix $n>2$ and let 
$\pi_{n+1}$ be a partition having 
a singleton block $\{n+1\}$. There is a chance that some $2\leq r\leq n$ balls are sampled from $\pi_{n+1}$ and added in the box $\{n+1\}$.
In this case the restriction of ${\rm SA}_{n+1}$ to $[n]$ {\it creates} a novel nonsingleton box, which is not a legitime option
for ${\rm SA}_n$.

\par In \cite{RPS} it was shown that
 the unique stationary $[n]$-partition is the one given by  (\ref{recregen1}).

\par 
{\bf Example.}
When
$$q(n:1)={n\rho\over 1+n\rho}\,\,,~~~~q(n:n)={1\over 1+n\rho},$$
the operation will create a new singleton block with probability $q(n:1)$, and merge everything in one block with probability $q(n:n)$. 
So the stationary distributions will concentrate on hook partitions.
The decrement matrix for this chain is the same as for ${\rm SA}_n$.

\par 
{\bf Example.} When
\eq
q(n:m)={n\choose m}\frac{[\theta]_{n-m}m!}{[\theta+1]_{n-1}n},
\en
with $\theta=2\rho$, 
the invariant partition is  Ewens' with parameter $\theta$. 
The decrement matrix for this chain is different from the one for 
${\rm SA}_n$, which also leads to Ewens' distribution.

\subsection{Comparing decrement matrices}

In \cite{RCS} we found very similar recursions for entries of decrement matrix which characterizes a regenerative composition structure, hence a
regenerative partition structure in \cite{RPS}. According to \cite[Proposition 3.3]{RCS}, a non-negative matrix $q$ is the decrement matrix of 
some regenerative composition structure if and only if $q(1:1)=1$ and (\ref{recq}) holds for $1\le k\le b$. 
Comparing with Lemma \ref{l5} above, the difference from our recursions here is that we have a separate recursion for $q(b:1)$, 
and we have an extra term
$$
\frac{2}{b+1}q(b+1:2)q(b:k)
$$
in right hand side of recursions for $q(b:k)$, $k\ge 2$. Both of them are backward recursions. 
For the purpose of illustration, 
suppose we are given $q(4:k)$, $k=1,2,3,4$, the entries $q(b:\cdot)$ with $b\le 3$ of decrement matrix for regenerative composition structure would be:
\begin{align*}
q(3:3)&=\frac{4q(4:4)+q(4:3)}{4-q(4:1)}\,,\\
q(3:2)&=\frac{3q(4:3)+2q(4:2)}{4-q(4:1)}\,,\\
q(3:1)&=\frac{2q(4:2)+3q(4:1)}{4-q(4:1)}\,,
\end{align*}
\begin{align*}
q(2:2)&=\frac{3q(3:3)+q(3:2)}{3-q(3:1)}=\frac{6q(4:4)+3q(4:3)+q(4:2)}{6-3q(4:1)-q(4:2)}\,,\\
q(2:1)&=\frac{2q(3:2)+2q(3:1)}{3-q(3:1)}=\frac{3q(4:3)+4q(4:2)+3q(4:1)}{6-3q(4:1)-q(4:2)}\,.
\end{align*}
While for decrement of the partition structure studied here, we have
\begin{align*}
q(3:3)&=\frac{4q(4:4)+q(4:3)}{4-q(4:1)-2q(4:2)}\,,\\
q(3:2)&=\frac{3q(4:3)+2q(4:2)}{4-q(4:1)-2q(4:2)}\,,\\
q(3:1)&=\frac{3q(4:1)}{4-q(4:1)-2q(4:2)},
\end{align*}
\begin{align*}
q(2:2)&=\frac{3q(3:3)+q(3:2)}{3-q(3:1)-2q(3:2)}=\frac{6q(4:4)+3q(4:3)+q(4:2)}{6-3q(4:1)-5q(4:2)-3q(4:3)}\,,\\
q(2:1)&=\frac{2q(3:1)}{3-q(3:1)-2q(3:2)}=\frac{3q(4:1)}{6-3q(4:1)-5q(4:2)-3q(4:3)}\,.
\end{align*}

\section{Comparison with Markovian fragmentations}
\label{9}

The theory of homogenous and self-similar Markovian fragmentation 
processes due to Bertoin \cite{BertoinBook} is formulated much
like the present theory of coalescents in terms of consistent
partition-valued processes.
Ford \cite[Proposition 41]{ford} provides a sampling consistency condition
for decrement matrices associated with discrete fragmentation processes which 
is an extremely close relative of our Lemma \ref{l5}.  
The article \cite{HMPW} provides an integral representation for such
decrement matrices, analogous to our results for the decrement matrices
associated with regenerative partition structures and with 
Markovian coalescents, and embeds Ford's result in the broader context 
of continuous time fragmentation processes and continuum random trees.
A missing element of the fragmentation discussion is some way of
deriving a partition structure by a recursion like 
\re{rec1} or \re{recregen}. But we expect such a partition structure and
an associated recursion may be associated with a suitably defined
Markovian fragmentation with freeze, such as that introduced in
\cite{GY}.


\begin{thebibliography}{10}

\bibitem{aldous85}
D.~J. Aldous.
\newblock Exchangeability and related topics.
\newblock In {\em \'Ecole d'\'et\'e de probabilit\'es de Saint-Flour,
  XIII---1983}, {\em Lecture Notes Math.} 1117: 1--198.
  Springer, Berlin, 1985.

\bibitem{BertoinBook} 
J. Bertoin, 
\newblock Random fragmentation and coagulation processes.
\newblock Cambridge University Press, to appear.

\bibitem{bertoingold04}
J. Bertoin and C. Goldschmidt.
\newblock Dual random fragmentation and coagulation and an application to the
  genealogy of {Y}ule processes.
\newblock In {\em Mathematics and computer science. III}, Trends Math., 
  295--308. Birkh\"auser, Basel, 2004.

\bibitem{boltszni98}
E. Bolthausen and A.~S. Sznitman.
\newblock On {R}uelle's probability cascades and an abstract cavity method.
\newblock {\em Comm. Math. Phys.} 197(2): 247--276, 1998.

\bibitem{canning74}
C. Cannings.
\newblock The latent roots of certain {M}arkov chains arising in genetics: a
  new approach. {I}. {H}aploid models.
\newblock {\em Advances in Appl. Probability} 6: 260--290, 1974.

\bibitem{donggoldmartin05}
R. Dong, C. Goldschmidt, and J.~B. Martin.
\newblock {Coagulation-fragmentation duality, Poisson-Dirichlet distributions
  and random recursive trees}.
\newblock Preprint (2005), arXiv:math.PR/0507591


\bibitem{donnellyjoyce91}
P. Donnelly and P. Joyce.
\newblock Consistent ordered sampling distributions: characterization and
  convergence.
\newblock {\em Adv. in Appl. Probab.} 23(2): 229--258, 1991.

\bibitem{donnellytavare86}
P. Donnelly and S. Tavar{\'e}.
\newblock The ages of alleles and a coalescent.
\newblock {\em Adv.  Appl. Probab.} 18(1): 1--19, 1986.

\bibitem{evanspitman98}
S.~N. Evans and J. Pitman.
\newblock Construction of {M}arkovian coalescents.
\newblock {\em Ann. Inst. H. Poincar\'e Probab. Statist.} 34(3): 339--383,
  1998.

\bibitem{ewens72}
W.~J. Ewens.
\newblock The sampling theory of selectively neutral alleles.
\newblock {\em Theoret. Population Biology}, 3: 87--112; erratum, ibid. 3
  (1972), 240; erratum, ibid. 3 (1972), 376, 1972.

\bibitem{ford}
D.~J. Ford.
\newblock Probabilities on cladograms: introduction to the alpha model.
\newblock Preprint (2005), arXiv:math.PR/0511246

\bibitem{gnedin97}
A.~V. Gnedin.
\newblock The representation of composition structures.
\newblock {\em Ann. Probab.} 25(3): 1437--1450, 1997.

\bibitem{selfsim} A.~Gnedin and J.~Pitman. Markov and self-similar composition structures,
{\it Zapiski Pomi} 326: 59--84. Available at
{\tt http://www.pdmi.ras.ru/znsl/2005/v326.html}

\bibitem{RCS}
A. Gnedin and J. Pitman.
\newblock Regenerative composition structures.
\newblock {\em Ann. Probab.} 33(2): 445--479, 2005.

\bibitem{RPS}
A. Gnedin and J. Pitman.
\newblock Regenerative partition structures.
\newblock {\em Electron. J. Combin.}, 11(2): Research Paper 12, 21 pp.
  (electronic), 2004/05.

\bibitem{gnedinpitman06}
A. Gnedin and J. Pitman.
\newblock {Moments of convex distribution functions and completely alternating
  sequences}, 2006,
arXiv:math.PR/0602091.


\bibitem{GY} 
A. Gnedin and Y. Yakubovich. 
\newblock Recursive Partition Structures,
2005, arXiv:math.PR/0510305.

\bibitem{HMPW} 
B. Haas, G. Miermont, J. Pitman and M. Winkel
\newblock Asymptotics of discrete fragmentation trees and applications to phylogenetic models, 2006,
\newblock in preparation.

\bibitem{kingman78a}
J.~F.~C. Kingman.
\newblock The representation of partition structures.
\newblock {\em J. London Math. Soc. (2)} 18(2): 374--380, 1978.

\bibitem{kingman82a}
J.~F.~C. Kingman.
\newblock The coalescent.
\newblock {\em Stochastic Process. Appl.} 13(3):235--248, 1982.

\bibitem{kingman82c}
J.~F.~C. Kingman.
\newblock Exchangeability and the evolution of large populations.
\newblock In {\em Exchangeability in probability and statistics (Rome, 1981)},
  97--112, North-Holland, Amsterdam, 1982.

\bibitem{kingman82b}
J.~F.~C. Kingman.
\newblock On the genealogy of large populations.
\newblock {\em J. Appl. Probab.} (Special Vol. 19A): 27--43, 1982.
\newblock Essays in statistical science.

\bibitem{moehle1}
M. M{\"o}hle.
\newblock On sampling distributions for coalescent processes with simultaneous
  multiple collisions.
\newblock {\em  Bernoulli} 12(2): 35--53, 2006.

\bibitem{moehle2}
M. M{\"o}hle.
\newblock On a class of non-regenerative sampling distributions.
\newblock {\em Combinatorics, Probability and Computing} 2005, to appear.

\bibitem{moehlesagitove01}
M. M{\"o}hle and S. Sagitov.
\newblock A classification of coalescent processes for haploid exchangeable
  population models.
\newblock {\em Ann. Probab.} 29(4): 1547--1562, 2001.

\bibitem{moran58}
P.~A.~P. Moran.
\newblock Random processes in genetics.
\newblock {\em Proc. Camb. Phil. Soc.} 54: 60--71, 1958.

\bibitem{Nordborg} 
M. Nordborg. 
\newblock Coalescent theory. 
\newblock In D.J. Balding et al (eds) {\it Handbook
of statistical genetics}, 179-208, Wiley, NY, 2001.

\bibitem{CSP}
J. Pitman.
\newblock {\em Combinatorial stochastic processes}.
\newblock  (Lecture notes for St. Flour course, July 2002),
{\it Springer L. Notes Math.}, 2006, to appear.
Available via {\tt http://www.stat.Berkeley.edu}

\bibitem{pitman95}
J. Pitman.
\newblock Exchangeable and partially exchangeable random partitions.
\newblock {\em Probab. Theory Related Fields} 102(2): 145--158, 1995.

\bibitem{pitman97}
J. Pitman.
\newblock Partition structures derived from {B}rownian motion and stable
  subordinators.
\newblock {\em Bernoulli} 3(1): 79--96, 1997.

\bibitem{lambda}
J. Pitman.
\newblock Coalescents with multiple collisions.
\newblock {\em Ann. Probab.} 27(4): 1870--1902, 1999.

\bibitem{pitmanyor97}
J. Pitman and M. Yor.
\newblock The two-parameter {P}oisson-{D}irichlet distribution derived from a
  stable subordinator.
\newblock {\em Ann. Probab.} 25(2):855--900, 1997.

\bibitem{sagitov}
S. Sagitov.
\newblock The general coalescent with asynchronous mergers of ancestral lines.
\newblock {\em J. Appl. Probab.} 36(4): 1116--1125, 1999.

\bibitem{sagitov03}
S. Sagitov.
\newblock Convergence to the coalescent with simultaneous multiple mergers.
\newblock {\em J. Appl. Probab.} 40(4): 839--854, 2003.

\bibitem{jason00a}
J. Schweinsberg.
\newblock A necessary and sufficient condition for the {$\Lambda$}-coalescent
  to come down from infinity.
\newblock {\em Electron. Comm. Probab.} 5:1--11 (electronic), 2000.

\bibitem{jason00}
J. Schweinsberg.
\newblock Coalescents with simultaneous multiple collisions.
\newblock {\em Electron. J. Probab.} 5: Paper no.\ 12, 50 pp. (electronic),
  2000.
  
 \bibitem{tavare84}
S. Tavar{\'e}.
\newblock Line-of-descent and genealogical processes, and their applications in
  population genetics models.
\newblock {\em Theoret. Population Biol.}  26(2): 119--164, 1984.

\bibitem{watterson76}
G.~A. Watterson.
\newblock Reversibility and the age of an allele.
\newblock {\em Theoret. Population Biology} 10: 239--253, 1976.

\bibitem{watter84}
G.~A. Watterson.
\newblock Lines of descent and the coalescent.
\newblock {\em Theoret. Population Biol.} 26(1): 77--92, 1984.

\bibitem{young95}
J.~E. Young.
\newblock Partition-valued stochastic processes with applications.
\newblock {\em U.C. Berkeley Ph.D. Thesis}, 1995.

\end{thebibliography}
\end{document}